\documentclass[a4paper,twoside,10pt]{article}
\usepackage[a4paper,left=3cm,right=3cm, top=3cm, bottom=3cm]{geometry}
\usepackage[T1]{fontenc}
\usepackage{soul}
\usepackage{mathrsfs}
\usepackage{amsmath}
\usepackage{amsthm}
\usepackage{amssymb}
\usepackage{hyperref}
\usepackage{cite}
\usepackage{color}
\usepackage{doi}
\usepackage{cancel}
\usepackage{appendix}
\usepackage{tikz}
\usepackage{soul}
\theoremstyle{plain}
\newtheorem{theorem}{Theorem}[section]
\newtheorem{corollary}[theorem]{Corollary}
\newtheorem{lemma}[theorem]{Lemma}

\theoremstyle{remark}
\newtheorem{remark}{Remark}

\newcommand{\eremk}{\hbox{}\hfill\rule{0.8ex}{0.8ex}}

\DeclareTextCompositeCommand{\u}{T1}{i}{\u\imath}

\let\div\relax
\DeclareMathOperator{\div}{\nabla\cdot}

\DeclareMathOperator{\divbfh}{\boldsymbol \nabla_{\h} \cdot}

\DeclareMathOperator{\nablabold}{\boldsymbol \nabla}
\newcommand{\nablaboldh}{\nablabold_{\h}}
\DeclareMathOperator{\curlbold}{\nablabold\times}
\DeclareMathOperator{\Dtwo}{\underline{\bf D}^2}
\DeclareMathOperator{\Dtwoh}{\underline{\bf D}^2_{\h}}
\DeclareMathOperator{\Laplacetwo}{\Delta^2}

\newcommand{\Deltah}{\Delta_{\h}}

\newcommand{\Norm}[1]{{\left\|{#1} \right\|}}
\newcommand{\SemiNorm}[1]{{\left|{#1} \right|}}
\newcommand{\jump}[1]{\left[\!\left[#1\right]\!\right]}
\newcommand{\avg}[1]{\{#1\}}
\newcommand{\Normth}[1]{{\left\vert\kern-0.25ex\left\vert\kern-0.25ex\left\vert #1 \right\vert\kern-0.25ex\right\vert\kern-0.25ex\right\vert}}

\newcommand{\Nbb}{\mathbb N}
\newcommand{\Pbb}{\mathbb P}
\newcommand{\Rbb}{\mathbb R}

\let\H\relax
\DeclareMathOperator{\H}{H}
\let\L\relax
\DeclareMathOperator{\L}{L}
\let\W\relax
\DeclareMathOperator{\W}{W}
\newcommand{\Ccal}{\mathcal C}

\DeclareMathOperator{\Ccalz}{\Ccal^0}
\DeclareMathOperator{\Mcal}{\mathcal M}

\newcommand{\Lbf}{\mathbf L}
\newcommand{\Lbb}{\underline{\mathbf L}}
\newcommand{\Wbf}{\mathbf W}
\newcommand{\Wbb}{\underline{\mathbf W}}

\DeclareMathOperator{\VDG}{V_{\DG}}
\DeclareMathOperator{\VLL}{V_{\Delta^2}}
\DeclareMathOperator{\VCR}{V_{\CR}}
\DeclareMathOperator{\VCcalz}{V_{\Ccalz}}
\DeclareMathOperator{\VMcal}{V_{\Mcal}}

\DeclareMathOperator{\h}{h}
\DeclareMathOperator{\E}{K}
\DeclareMathOperator{\F}{F}
\DeclareMathOperator{\T}{T}
\newcommand{\e}{e}
\DeclareMathOperator{\D}{D}
\newcommand{\hE}{\h_{\E}}

\newcommand{\hF}{\h_{\F}}
\newcommand{\htildeF}{\widetilde{\h}_{\F}}

\newcommand{\hD}{\h_{\D}}
\newcommand{\hOmega}{\h_\Omega}

\newcommand{\FcalE}{\mathcal F^{\E}}

\DeclareMathOperator{\Fcaln}{\mathcal F_n}
\DeclareMathOperator{\FcalnI}{\mathcal F_n^I}
\DeclareMathOperator{\FcalnID}{\mathcal F_n^{ID}}
\DeclareMathOperator{\FcalnIN}{\mathcal F_n^{IN}}
\DeclareMathOperator{\FcalnB}{\mathcal F_n^B}
\DeclareMathOperator{\FcalnD}{\mathcal F_n^D}
\DeclareMathOperator{\FcalnN}{\mathcal F_n^N}

\newcommand{\taun}{\mathcal T_n}

\newcommand{\tautildeEF}{\widetilde{\mathcal T}^{\EF}}

\newcommand{\GammaN}{\Gamma_N}
\newcommand{\GammaD}{\Gamma_D}
\newcommand{\Gammatilde}{\widetilde\Gamma}

\newcommand{\FET}{\F_{\E}^{\T}}
\newcommand{\TfrakE}{\mathfrak T_{\E}}
\newcommand{\TfrakEF}{\mathfrak T_{\EF}}

\newcommand{\EF}{\E_{\F}}
\newcommand{\TEF}{\T_{\EF}}
\newcommand{\FTEF}{\F_{\TEF}}
\newcommand{\EfrakEF}{\mathfrak E _{\EF}}
\newcommand{\EfrakTEF}{\mathfrak E _{\TEF}}

\newcommand{\EoF}{\E_{1,\F}}
\newcommand{\EtF}{\E_{2,\F}}

\newcommand{\nbf}{\mathbf n}
\newcommand{\nbfE}{\nbf_{\E}}
\newcommand{\nbfOmega}{\nbf_\Omega}
\DeclareMathOperator{\nbfF}{\mathbf n_{\F}}
\DeclareMathOperator{\nbfFcaln}{\mathbf n_{\Fcaln}}

\newcommand{\nbfEoF}{\nbf_{\EoF}}
\newcommand{\nbfEtF}{\nbf_{\EtF}}
\newcommand{\nbfEF}{\nbf_{\EF}}

\newcommand{\vbf}{\mathbf v}

\DeclareMathOperator{\zerobf}{\bf 0}

\newcommand{\vh}{v_{\h}}
\DeclareMathOperator{\uDG}{u_{\DG}}
\DeclareMathOperator{\vDG}{v_{\DG}}
\DeclareMathOperator{\uCR}{u_{\CR}}
\DeclareMathOperator{\vCR}{v_{\CR}}
\DeclareMathOperator{\uCcalz}{u_{\Ccalz}}
\DeclareMathOperator{\vCcalz}{v_{\Ccalz}}
\DeclareMathOperator{\uMcal}{u_{\Mcal}}
\DeclareMathOperator{\vMcal}{v_{\Mcal}}

\newcommand{\sigmabold}{\boldsymbol \sigma}

\newcommand{\xibold}{\boldsymbol{\xi}}

\newcommand{\Psibold}{\boldsymbol \Psi}

\newcommand{\PizE}{\Pi^{0,\E}_0}
\newcommand{\PizF}{\Pi^{0,\F}_0}
\newcommand{\PizFcaln}{\Pi^{0}_{\Fcaln}}
\newcommand{\Piztaun}{\Pi^{0}_{\taun}}

\DeclareMathOperator{\p}{p}
\DeclareMathOperator{\q}{q}
\let\k\relax
\DeclareMathOperator{\k}{k}
\DeclareMathOperator{\s}{s}

\newcommand{\ostar}{1^*}
\newcommand{\pprime}{\p'}
\newcommand{\psharp}{\p^\sharp}
\newcommand{\pstar}{\p^*}
\newcommand{\postar}{\p \ostar}
\newcommand{\pprimeostar}{\pprime\ostar}
\newcommand{\psharpprime}{(\p^\sharp)'}
\newcommand{\postarprime}{(\p \ostar)'}
\newcommand{\postarostar}{\postar\ostar}
\newcommand{\postarostarprime}{(\postar\ostar)'}
\newcommand{\pstarprime}{(\p^*)'}
\newcommand{\pstarsharp}{(\p^*)^\sharp}
\newcommand{\pstarstar}{(\p^*)^*}
\newcommand{\pstarostar}{\pstar\ostar}
\newcommand{\pstarostarprime}{(\pstar\ostar)'}
\newcommand{\pstarprimesharp}{((\p^*)')^\sharp}
\newcommand{\pstarprimesharpprime}{(((\p^*)')^\sharp)'}

\newcommand{\pstarprimestar}{((\p^*)')^*}

\newcommand{\qprime}{\q'}

\DeclareMathOperator{\sprime}{\s'}
\DeclareMathOperator{\ssharp}{\s^\sharp}
\DeclareMathOperator{\sstar}{\s^*}

\DeclareMathOperator{\DG}{DG}
\DeclareMathOperator{\CR}{CR}

\DeclareMathOperator{\ap}{a_{\p}}

\DeclareMathOperator{\aCR}{a_{\CR}}
\DeclareMathOperator{\apCR}{a_{\p,\CR}}
\DeclareMathOperator{\btwo}{b_2}

\DeclareMathOperator{\btwoDG}{b_{2,\DG}}

\DeclareMathOperator{\btildetwo}{\widetilde b_2}

\DeclareMathOperator{\btildetwoDG}{\widetilde b_{2,\DG}}

\DeclareMathOperator{\btwoCcalz}{b_{2,\Ccalz}}

\DeclareMathOperator{\btwoMcal}{b_{2,\Mcal}}

\DeclareMathOperator{\cp}{c_{\p}}
\DeclareMathOperator{\cpMcal}{c_{\p,\Mcal}}

\DeclareMathOperator{\qkmoF}{q_{\k-1}^{\F}}

\DeclareMathOperator{\CHo}{C_{H,1}}
\DeclareMathOperator{\CHt}{C_{H,2}}
\DeclareMathOperator{\CBcal}{C_{\mathcal B}}
\DeclareMathOperator{\CPS}{C_{PS}}
\DeclareMathOperator{\CSob}{C_{Sob}}
\DeclareMathOperator{\CMT}{C_{MT}}
\DeclareMathOperator{\CTR}{C_{TR}}
\DeclareMathOperator{\CtildeTR}{\widetilde C_{TR}}
\DeclareMathOperator{\CBMoo}{C_{BM,1,1}}
\DeclareMathOperator{\CBMot}{C_{BM,1,2}}
\DeclareMathOperator{\CBMto}{C_{BM,2,1}}
\DeclareMathOperator{\CBMtt}{C_{BM,2,2}}
\DeclareMathOperator{\CBMtth}{C_{BM,2,3}}
\DeclareMathOperator{\CBMtf}{C_{BM,2,4}}
\DeclareMathOperator{\CtildeBMoo}{\widetilde C_{BM,1,1}}
\DeclareMathOperator{\CtildeBMot}{\widetilde C_{BM,1,2}}
\DeclareMathOperator{\CtildeBMoth}{\widetilde C_{BM,1,3}}
\DeclareMathOperator{\CtildeBMto}{\widetilde C_{BM,2,1}}
\DeclareMathOperator{\CtildeBMtt}{\widetilde C_{BM,2,2}}
\DeclareMathOperator{\CtildeBMtth}{\widetilde C_{BM,2,3}}
\DeclareMathOperator{\CtildeBMtho}{\widetilde C_{BM,3,1}}
\DeclareMathOperator{\CtildeBMtht}{\widetilde C_{BM,3,2}}
\DeclareMathOperator{\CtildeBMthth}{\widetilde C_{BM,3,3}}
\DeclareMathOperator{\CtildeBMfo}{\widetilde C_{BM,4,1}}
\DeclareMathOperator{\CtildeBMft}{\widetilde C_{BM,4,2}}
\DeclareMathOperator{\CtildeBMfth}{\widetilde C_{BM,4,3}}
\DeclareMathOperator{\CEG}{C_{EG}}
\DeclareMathOperator{\Csigma}{C_{\sigma}}
\DeclareMathOperator{\Ctau}{C_{\tau}}

\title{\small Sobolev--Poincar\'e inequalities for piecewise $\W^{2,\p}$ functions over general polytopic meshes}
\author{\footnotesize{Michele Botti\thanks{MOX, Department of Mathematics, Politecnico di Milano, 20133 Milano, Italy (michele.botti@polimi.it)},   
Lorenzo Mascotto\thanks{Department of Mathematics and Applications, University of Milano-Bicocca, 20125 Milan, Italy (lorenzo.mascotto@unimib.it); 
Faculty of Mathematics, University of Vienna, 1090 Vienna, Austria;
IMATI-CNR, 27100 Pavia, Italy}}}
\date{}

\begin{document}
\maketitle

\begin{abstract}
\noindent We establish Sobolev--Poincar\'e inequalities
for piecewise $\W^{2,\p}$ functions in two and three dimensions
over families of general polytopic
(including shape-regular simplicial and Cartesian) meshes,
allowing for degenerate facets.
The proposed inequalities feature control by either the broken Hessian or,
more naturally for general fourth-order problems, the broken Laplacian.
Important tools in their derivation are
a Helmholtz decomposition in Banach--Sobolev spaces
and Sobolev--Poincar\'e inequalities
for piecewise $\W^{1,\p}$ functions;
we avoid the use of any inverse inequality.
We provide estimates with fairly general Lebesgue exponents.
Applications to the stability analysis
of nonconforming discretizations for
second- and fourth-order, linear and nonlinear problems are exhibited.


\end{abstract}

\section{Introduction}

\paragraph*{State-of-the-art.}
Poincar\'e and Sobolev--Poincar\'e inequalities for functions
in broken Sobolev spaces are a crucial tool
in the analysis of nonconforming methods.
For instance, they are used to prove a priori estimates
for the discrete solution in terms of the data
\cite[Remark 1.15]{Botti-Mascotto:2026};
to derive upper bounds on the condition number
for the resulting nonconforming finite element
and preconditioned systems
\cite[Lemma 3.2]{Brenner-Owens:2007},
\cite[Theorem 5.2]{Carrero-Cockburn-Schoetzau:2006}; 
to establish error estimates in weaker norms
via the Aubin--Nitsche duality argument
for nonlinear problems
\cite[Theorem 4.10]{Gudi-Pani:2007}.

Such inequalities have been thoroughly studied
at least since the 1970s~\cite{Temam:2001};
see, among other references,
\cite{Arnold:1982, Brenner-Scott:2008, Brenner:2003, Vohralik:2005, Lasis-Suli:2003}. 
Most of the literature is concerned
with broken (with respect to sequences of ``regular'' meshes)
$\H^1$ and~$\W^{1,\p}$ spaces, $\p\in[1,\infty)$.
Adopting the notation of Section~\ref{section:assumption-standard-results} below,
such inequalities on a domain~$D$ take the form:
there exists a positive constant $C$ such that,
given~$\h$ the mesh-size function of the given mesh~$\taun$
and for all $v$ in $\W^{1,\p}(\taun)$,
\begin{equation} \label{eq:abstract-W1p-inequalities}
\Norm{v}_{\L^{\p}(D)}
\le C \Big( \Norm{\nablaboldh v}_{\Lbf^{\p}(D)} 
    + \Norm{\h^{\frac1{\p}-1}\text{jump}(v)}
        _{\L^{\p}(\text{skeleton})} \Big).
\end{equation}
Above, ``jump$(v)$'' represents either the jump of~$v$
or its projection onto constants over the mesh skeleton;
this interface term is added to the broken gradient seminorm
to recover a norm on the right-hand side
of~\eqref{eq:abstract-W1p-inequalities}.
These norms are broken with respect
to sequences of shape-regular simplicial meshes
or ``regular'' polytopic meshes
(e.g., in the sense of \cite[Section~4]{Brenner-Wang-Zhao:2004}).
Usually, no small facets are considered
and  the constant~$C$ depends on a mesh regularity parameter.
Most available references establish (Sobolev-)Poincar\'e inequalities
for broken finite-dimensional spaces
(e.g., piecewise polynomials)
by employing inverse inequalities; as a result,
$C$ is sometimes assumed to depend
on the underlying finite dimensional space as well.
This is however not always the case:
\cite{Vohralik:2005} hinges purely upon direct estimates, 
while~\cite{Brenner:2003} uses inverse estimates
only on fixed, low-degree piecewise polynomials.

Our recent work~\cite{Botti-Mascotto:2026} advanced
the state-of-the-art in several directions.
First, whereas previous Lebesgue indices in~\eqref{eq:abstract-W1p-inequalities}
were nonsharp compared to continuous Sobolev embeddings,
the new inequalities in~\cite{Botti-Mascotto:2026}
achieve sharp exponents.
Second, we accounted for mixed boundary conditions
and considered mesh sequences featuring general element geometries,
including elements with degenerate facets.
Finally, these inequalities were established
in arbitrary spatial dimensions for domains with general topology.
The proofs rely exclusively on direct estimates
and the construction of regularized right-inverses
of the divergence operator under mixed boundary conditions.

The literature concerning higher-order Sobolev--Poincar\'e inequalities
for broken  spaces is more limited and most
existing results rely on inverse estimates, which we aim to avoid.
One notable reference that primarily employs direct estimates
is~\cite{Brenner-Wang-Zhao:2004}, where inverse estimates are used
only for polynomial spaces of fixed, low degree.
In the Hilbert case (the Banach setting being analogous)
and on ``regular'' meshes, the strategy
in~\cite{Brenner-Wang-Zhao:2004} takes the following form:
there exists a positive constant $C$ such that,
given~$\h$ the mesh-size function of the given mesh~$\taun$
and for all $v$ in $\H^{2}(\taun)$,
\[
\begin{split}
\Norm{\nablaboldh v}_{\Lbf^2(\Omega)}
&   \le \Norm{\nablaboldh (v-v_{1})}_{\Lbf^2(\Omega)} + \Norm{\nablaboldh v_{1}}_{\Lbf^2(\Omega)}\\
&   \le C \big(\Norm{\h \Dtwoh v}_{\Lbb^2(\Omega)} + \Norm{\nablaboldh v_{1}}_{\Lbf^2(\Omega)}\big) \\
&   \le C \Big(\Norm{\Dtwoh v}_{\Lbb^2(\Omega)} 
    + \Norm{\h^{-\frac32} \text{jump}(v)}_{\L^2(\text{skeleton})}
    + \Norm{\h^{-\frac12} \text{jump}(\text{flux} \ v)}_{\Lbf^2(\text{skeleton})} \Big).
\end{split}
\]
Above, $v_1$ denotes, for instance, the piecewise linear
Lagrange interpolant of~$v$
and $\Dtwoh$ is the piecewise Hessian operator.
The second inequality relies on standard local interpolation estimates,
whereas the third inequality involves inverse estimates,
cf. \cite[Section~2]{Brenner-Wang-Zhao:2004},
based on the  approximation properties of enriching operators
onto HCT spaces over a suitable subtriangulation.

This approach has a major advantage: it extends
fairly easily to several other types of inequalities,
including broken $\H^1$~\cite{Brenner:2003}
and broken Korn~\cite{Brenner:2004} inequalities.
However, to the best of our knowledge,
a key limitation is that the $\L^2$ norm of the Hessian
appears on the right-hand side,
which is stronger than the $\L^2$ norm of the Laplacian.
The latter is more natural for analyzing nonconforming methods
applied to problems such as the biharmonic problem
with mixed boundary conditions~\cite{Dang:2006}
or the $\p$-biharmonic problem~\cite{Pryer:2014},
which arises in modeling nonlinear plate and beam bending,
viscoelasticity, and phase-field systems.

\paragraph*{Main results and goals of the paper.}
The main goal of this paper is to establish Sobolev--Poincar\'e inequalities
for piecewise $\W^{2,\p}$ functions in two and three dimensions,
advancing the current state-of-the-art in several key directions.
\begin{itemize}
\item The inequalities hold true
on sequences of fairly general meshes,
including those with degenerate facets,
thereby accommodating meshes produced by coarsening steps
in adaptive mesh refinement algorithms.
\item They accommodate different
Lebesgue exponents to broaden applicability,
enabling, for example, stability analysis
of nonconforming methods with reduced data regularity
for linear and nonlinear fourth-order problems.
\item They apply directly to functions
in broken Sobolev spaces rather than
being restricted to discrete finite element subspaces,
thereby avoiding inverse estimates altogether.
\item Where possible, they feature the $\L^{\p}$
norm of the Laplacian
rather than the Hessian, providing a more natural setting
for analyzing general fourth-order problems.
\end{itemize}

We provide different inequalities,
each with its own advantages and limitations;
several applications to the stability analysis
of nonconforming methods for nonlinear problems are given.
Our strategy relies on a combination of tools,
one of the most important being
a Helmholtz decomposition of the broken gradient
of piecewise $\W^{2,\p}$ functions.

\paragraph*{Outline of the paper and list of relevant results.}
The paper is structured as follows.
Section~\ref{section:assumption-standard-results}
presents basic notation, the definition of ``regular'' mesh sequences,
auxiliary inequalities, and Sobolev--Poincar\'e inequalities
for piecewise $\W^{2,\p}$ functions involving
the $\L^{\p}$ norm of the broken Hessian,
alongside applications to the stability
of nonconforming methods for the $\p$-Laplace problem
and the biharmonic problem in Hessian--Hessian formulation.
Section~\ref{section:main} derives Sobolev--Poincar\'e inequalities
featuring the $\L^{\p}$ norm of the broken Laplacian,
introduces a key Helmholtz decomposition,
and addresses applications to the $\p$-biharmonic problem
and the biharmonic problem in Laplace--Laplace formulation
with mixed boundary conditions.
Finally, Appendix~\ref{appendix:regularity-Poisson} collects regularity results
in Banach-Sobolev spaces for two- and three-dimensional Poisson problems,
which dictate the regularity of the vector fields in the Helmholtz decomposition.

For the readers' ease, we list here what we believe
the relevant results of this work are.
\begin{itemize}
\item Theorem~\ref{theorem:BM-W1p} recalls the Sobolev--Poincar\'e inequalities
for piecewise $\W^{1,\p}$ functions from~\cite{Botti-Mascotto:2026}.
\item Corollary~\ref{corollary:trivial-W2p-inequalities}
exhibits ``standard'' Sobolev--Poincar\'e inequalities
for piecewise $\W^{2,\p}$ functions involving
the $\L^{\p}$ norm of the piecewise Hessian,
thus recovering results similar to those
in~\cite{Brenner-Wang-Zhao:2004}.
\item Applications of Theorem~\ref{theorem:BM-W1p}
and Corollary~\ref{corollary:trivial-W2p-inequalities}
to the stability analysis of nonconforming methods
for the $\p$-Laplace problem and the biharmonic problem
in Hessian--Hessian formulation are given
in Section~\ref{subsection:application-1}.
\item Theorem~\ref{theorem:sharp-ineq} exhibits
Sobolev--Poincar\'e inequalities
for piecewise $\W^{2,\p}$ functions involving
the $\L^{\p}$ norm of the piecewise Laplacian
with sharp Lebesgue indices,
yet with certain limitations in the range of~$\p$.
\item Theorem~\ref{theorem:nonsharp-ineq} exhibits
Sobolev--Poincar\'e inequalities
for piecewise $\W^{2,\p}$ functions involving
the $\L^{\p}$ norm of the piecewise Laplacian
with nonsharp Lebesgue indices
with milder limitations in the range of~$\p$.
\item Applications of Theorems~\ref{theorem:sharp-ineq}
and~\ref{theorem:nonsharp-ineq}
to the stability analysis of nonconforming methods
for the $\p$-biharmonic problem and the biharmonic problem
in Laplace--Laplace formulation with mixed boundary conditions
are given in Section~\ref{subsection:application-2}.
\item A Helmholtz decomposition, which is instrumental in the
derivation of Theorems~\ref{theorem:sharp-ineq}
and~\ref{theorem:nonsharp-ineq}, is given
in Section \ref{subsection:Helmholtz-dec}.
\end{itemize}

\section{Assumptions and standard results} \label{section:assumption-standard-results}
For a generic Lipschitz domain~$\D$ with diameter~$\hD$,
standard scalar, vector, and tensor Lebesgue spaces
of order~$\p\in[1,\infty]$ are denoted
by~$\L^{\p}(\D)$, $\Lbf^{\p}(\D)$, and $\Lbb^{\p}(\D)$;
one can also consider Lebesgue spaces on
the boundary~$\partial \D$ of~$\D$
or part thereof. An additional subscript~$0$ means
that we consider scalar, vector, and tensor fields
with zero average over the domain.
We consider standard differential operators:
the gradient~$\nablabold$; the Hessian~$\Dtwo$;
the divergence~$\div$; the curl~$\curlbold$;
the Laplacian~$\Delta$;
the bilaplacian~$\Laplacetwo$;
given $k$ in $\Nbb$, the $k$-tensor $D^{k}$
containing all mixed derivatives up to order~$k$.

Given~$(k,s)$ in $\{0,1,2\}\times[1,\infty)$,
consider the seminorms and norms
\footnote{With an abuse of notation, here we use the subscript
$\L^s(\D)$ also for vector and tensor fields.}
\[
\SemiNorm{v}_{\W^{k,s}(\D)}^s
:= \Norm{D^k v}_{\L^{s}(\D)}^s,
\qquad\qquad\qquad
\Norm{v}_{\W^{k,s}(\D)}^s
:= \sum_{\ell=0}^k (\hD^{\ell-k} \SemiNorm{v}_{\W^{\ell,s}(\D)})^s.
\]
The Sobolev spaces $\W^{k,s}(\D)$,
$(k,s)\in \Nbb \times [1,\infty)$
are the spaces of scalar~$\L^s(\D)$
functions with finite norm $\Norm{\cdot}_{\W^{k,s}(\D)}$;
Sobolev spaces of noninteger order can be defined via interpolation theory.
Since a trace theorem holds true,
cf. \cite[Section~3]{Ern-Guermond:2021},
given $\Gammatilde \subseteq \partial\D$
with nonzero measure, we also write
\[
\W^{1,\s}_{\Gammatilde} (\D)
:= \{ v\in \W^{1,\s} (\D) \mid v_{|\Gammatilde} =0 \},
\qquad\qquad\qquad
\W^{1,\s}_{0} (\D)
:= \W^{1,\s}_{\partial\D} (\D) .
\]
We shall be also considering the spaces
\[
\W^{2,\s}_{\Gammatilde} (\D)
:= \{ v\in \W^{2,\s}\cap\W^{1,\s}_{\Gammatilde}(\D)
    \mid \nablabold v_{|\Gammatilde} = \zerobf \},
\qquad\qquad\qquad
\W^{2,\s}_{0} (\D)
:= \W^{2,\s}_{\partial\D} (\D) .
\]\normalsize
Vector and tensor Sobolev spaces are denoted
by~$\Wbf^{k,s}(\D)$ and~$\Wbb^{k,s}(\D)$,
and corresponding notation in case a zero trace is assigned
on (part of) $\partial\D$.

Given~$\p\in[1,\infty)$, we set
\[
    \pprime :=  \frac{\p}{\p-1};
    \qquad
    \pstar := 
    \begin{cases}
    \frac{\p d}{d-\p}   & \text{if } \p < d,  \\
    \infty              & \text{otherwise};
    \end{cases}
    \qquad
    \psharp := 
    \begin{cases}
    \frac{\p(d-1)}{d-\p}    & \text{if } \p < d, \\
    \infty                  & \text{otherwise}.
    \end{cases}
\]
Henceforth, $\Omega\subset \Rbb^d$, $d=2,3$,
is a bounded, polytopic, Lipschitz domain
with diameter~$\hOmega$. Depending on the occasion,
we may either require
\begin{subequations}
\begin{align} 
\text{$\Omega$ is star-shaped with respect to a ball;} \label{eq:assumption-Omega-1} \\
\text{$\Omega$ admits a shape-regular
decomposition into simplices;} \label{eq:assumption-Omega-2} \\
\text{$\Omega$ has
$(d-1)$ Betti number $\beta_{d-1}(\Omega)=0$
in $\Rbb^d$} \label{eq:assumption-Omega-3}.
\end{align}
\end{subequations}
Note that \eqref{eq:assumption-Omega-1} implies~\eqref{eq:assumption-Omega-3}.
The boundary~$\partial\Omega=:\Gamma$ of~$\Omega$
is split into~$\Gamma=\GammaN\cup\GammaD$;
we assume that~$\GammaD$
has nonzero $(d-1)$-dimensional Hausdorff measure,
that~$\GammaD\cap\GammaN=\emptyset$,
and that~$\GammaD$ is closed in the topology of~$\Gamma$.

The remainder of this section is organized as follows.
Sequences of ``regular'' polytopic meshes
along with broken spaces are given in Section~\ref{subsection:meshes-broken};
auxiliary, classical inequalities are collected
in Section~\ref{subsection:inequalities};
Sobolev--Poincar\'e inequalities for piecewise $\W^{1,\p}$ functions
from~\cite{Botti-Mascotto:2026}
and Sobolev--Poincar\'e inequalities
for piecewise $\W^{2,\p}$ functions involving
the $\L^{\p}$ norm of the piecewise Hessian
are discussed in Section~\ref{subsection:trivial-inequalities};
their applications to the stability analysis of nonconforming methods
for the $\p$-Laplace problem and the biharmonic problem
in Hessian--Hessian formulation are given
in Section~\ref{subsection:application-1}.

\subsection{Meshes and broken spaces} \label{subsection:meshes-broken}
Let~$\{\taun\}$ be a family of meshes,
where each~$\taun$ consists of a finite set of disjoint,
closed, polytopic elements satisfying
$\overline{\Omega}=\bigcup_{\E\in\taun}\E$.
Differential operators acting elementwise over~$\taun$
have a subscript~$\h$; as an example, $\nablaboldh$
represents the broken gradient operator on~$\taun$.
For any element~$\E \in \taun$, its boundary and diameter
are denoted by~$\partial\E$ and~$\hE$, respectively.

With each~$\taun$, we associate the set of facets~$\Fcaln$
covering the mesh skeleton
so that~$\bigcup_{\E\in\taun}\partial\E = \bigcup_{\F\in\Fcaln}\F$.
Each facet~$\F \in \Fcaln$ is a hyperplanar, closed,
connected subset of~$\overline{\Omega}$
having positive $(d-1)$-dimensional Hausdorff measure,
categorized as follows:
\begin{itemize}
\item an internal facet (interface) if there exist
two distinct elements~$\EoF, \EtF \in \taun$
such that~$\F\subseteq\partial\EoF\cap\partial\EtF$;
\item a boundary facet if there exists
an element~$\EF \in \taun$ such that
$\F\subseteq\partial\EF\cap\Gamma$.
\end{itemize}
The sets of internal and boundary facets
are denoted by~$\FcalnI$ and~$\FcalnB$, respectively.
The boundary skeleton is assumed to conform
to the partition of~$\Gamma$
into~$\GammaD$ and~$\GammaN$,
meaning~$\FcalnB = \FcalnD \cup \FcalnN$, where
\[
\FcalnD := \{ \F\in\FcalnB \mid \F \subseteq \GammaD \},
\qquad\qquad\qquad
\FcalnN := \{ \F\in\FcalnB \mid \F \subseteq \GammaN \}.
\]
We also write $\FcalnID:=\FcalnI\cup\FcalnD$
and $\FcalnIN:=\FcalnI\cup\FcalnN$.
For an element~$\E \in \taun$,
$\FcalE$ represents its collection of boundary facets
and~$\nbfE$ denotes the outward unit normal vector to~$\partial\E$.
To each internal facet~$\F \in \FcalnI$,
we fix~$\nbfF$ as one
of the two unit normal vectors to~$\F$;
to each boundary facet~$\F \in \FcalnB$
with $\F\subset\partial\E$, we assign $\nbfF=\nbfE$.
We further define the skeleton normal
unit vector~$\nbfFcaln$ as $\nbfFcaln{}_{|\F}:=\nbfF$.

The mesh family is assumed to satisfy
the following regularity criteria~\cite{Cangiani-Dong-Georgoulis-Houston:2017}:
\begin{itemize}
\item for every~$\E \in \taun$,
there exists a partition~$\TfrakE$ of~$\E$
into nonoverlapping $d$-dimensional simplices;
\item there exists a uniform constant~$\gamma > 0$ such that,
for all~$\E \in \taun$ and every simplex
$\T \in \TfrakE$ with~$\partial\T\cap\partial\E=\FET$
being a $(d-1)$-dimensional simplex,
\begin{equation} \label{eq:mesh-regularity}
\gamma \hE \le d \vert\T\vert \vert\FET\vert^{-1}.
\end{equation}
\end{itemize}
On occasions, we shall also require that
the partition~$\TfrakE$ is~$\gamma$ shape-regular:
\begin{equation} \label{eq:mesh-regularity-extra}
\gamma \h_{\T}^{d} \le \vert \T \vert
\qquad\qquad\qquad\forall\T\in\TfrakE.
\end{equation}
Define
\[
\htildeF
\quad \text{ either equal to } \hF
\text{ or equal to }\quad
\begin{cases}
    \min(\h_{\E_1} , \h_{\E_2})
        & \text{if } \F\in\FcalnI,\; \F\subset\partial\E_1\cap\partial\E_2,\\
    \h_{\E}                     
        & \text{if } \F\in\FcalnD,\; \F \in\FcalE.
\end{cases}         
\]
We furthermore let~$\h$ to denote
the piecewise constant function
such that $\h_{|\E}=\hE$ for all $\E\in\taun$
and $\h_{|\F}=\htildeF$ for all $\F\in\Fcaln$.
Broken Sobolev spaces associated with~$\taun$ are given by
\[
\W^{k,s}(\taun)
:= \left\{ v \in \L^{s}(\Omega) \;\middle|\;
    v_{|\E} \in \W^{k,s}(\E)
    \quad \forall \E \in \taun \right\},
\qquad k=1,2.
\]
For any scalar function~$v \in \W^{1,\p}(\taun)$
and facet~$\F \in \Fcaln$, the jump operator is defined as
(with analogous definitions
if vectors and tensors are considered)
\[
\jump{v}_{\F}
= \jump{v}
:= \begin{cases}
(v_{|\EoF} \nbfEoF \cdot \nbfF + v_{|\EtF} \nbfEtF \cdot \nbfF)_{|\F}
&   \text{if } \F\in\FcalnI,\ \F \subset\partial\EoF\cap\partial\EtF, \\
(v_{|\EF} \nbfEF \cdot \nbfF)_{|\F}
&   \text{if } \F\in\FcalnB,\ \F \subset \partial\EF\cap\Gamma,
\end{cases}
\]
while, for vector fields~$\vbf \in \Wbf^{1,\p}(\taun)$,
the normal and tangential jump operators take the form
(with analogous definitions if tensors are considered)
\[
\jump{\vbf}_{N,\F}
= \jump{\vbf}_{N}
:= \begin{cases}
(\nbfEoF\cdot\vbf_{|\EoF}  + \nbfEtF\cdot\vbf_{|\EtF})_{|\F}
&   \text{if } \F\in\FcalnI,\ \F \subset\partial\EoF\cap\partial\EtF, \\
(\nbfEF\cdot\vbf_{|\EF})_{|\F}
&   \text{if } \F\in\FcalnB,\ \F \subset \partial\EF\cap\Gamma.
\end{cases}
\]
and
\[
\jump{\vbf}_{T,\F}
= \jump{\vbf}_{T}
:= \begin{cases}
(\nbfEoF\times\vbf_{|\EoF}  + \nbfEtF\times\vbf_{|\EtF})_{|\F}
&   \text{if } \F\in\FcalnI,\ \F \subset\partial\EoF\cap\partial\EtF, \\
(\nbfEF\times\vbf_{|\EF})_{|\F}
&   \text{if } \F\in\FcalnB,\ \F \subset \partial\EF\cap\Gamma.
\end{cases}
\]
We also define average operators
(with analogous definitions for vectors and tensors)
\[
\avg{v}_{\F} = \avg{v} :=
\begin{cases}
\frac12 (v_{|\EoF}+ v_{|\EtF})_{|\F}
    &  \text{if } \F\in\FcalnI,\ \F \subset\partial\EoF\cap\partial\EtF,\\
(v_{|\EF})_{|\F}
    & \text{if } \F\in\FcalnB,\ \F \subset \partial\EF\cap\Gamma.
\end{cases}
\]
Finally, we define piecewise $\L^2$ projectors
onto constants $\PizF: \L^1(\F) \to \Pbb_{0}(\F)$
and $\PizE: \L^1(\E) \to \Pbb_0(\E)$.
The corresponding global projectors
over~$\Fcaln$ and~$\taun$ are~$\PizFcaln$ and~$\Piztaun$.
With an abuse of notation, we also write
\[
(u,v)_{\Fcaln}
:= \sum_{\F\in\Fcaln} (u,v)_{\L^2(\F)},
\qquad\qquad\qquad
\Norm{v}_{\Fcaln}:=(v,v)_{\Fcaln}.
\]
Analogous definition will be used for vector and tensors
and for subsets of~$\Fcaln$.

\subsection{Auxiliary inequalities} \label{subsection:inequalities}

The following result is standard;
cf., e.g., \cite[Chapter~3]{Ern-Guermond:2021}.
\begin{lemma}[Poincar\'e inequality]
Given a Lipschitz domain~$\D$ with diameter~$\hD$,
there exists a positive constant~$\CPS$ independent of~$\hD$
such that, for all $v$ in $\W^{1,\s}(\D)$, $\s\in[1,\infty)$,
with either zero trace on part of the boundary of~$\D$,
or zero average on part of~$\D$ or its boundary,
\begin{equation} \label{eq:standard-Poincare-Steklov}
    \Norm{v}_{\L^{\s}(\D)}
    \le \CPS(\s) \hD \SemiNorm{v}_{\W^{1,\s}(\D)}.
\end{equation}
\end{lemma}

Also the following result is fairly standard;
cf., e.g., \cite[Chapter~2]{Ern-Guermond:2021}:
\begin{lemma}[Gagliardo--Nirenberg--Sobolev embedding]
Given a Lipschitz domain~$\D$ with diameter~$\hD$,
there exists a positive constant~$\CSob$ independent of~$\hD$
such that, for~$\s$ in $[1,d)$,
\begin{equation} \label{eq:GNS-inequality}
    \Norm{v}_{\L^{\sstar}(\D)}
    \le \CSob(\s,\D) \Norm{v}_{\W^{1,\s}(\D)}
    \qquad\qquad\qquad\forall v \in \W^{1,\s}(\D).
\end{equation}
\end{lemma}

Next, we recall a local trace inequality;
cf., e.g., \cite[Corollary 1.3]{Botti-Mascotto:2026}.
\begin{lemma}[Sobolev trace inequality]
Let~$\E$ be an element in~$\taun$ and $s$ be in $[1,d)$.
Then, there exists a positive constant~$\CTR$
independent of~$\hE$ such that
\begin{equation} \label{eq:local-Sobolev-trace}
\Norm{v}_{\L^{\ssharp}(\partial\E)}
\le \CTR(\s,\gamma) \big( \Norm{v}_{\L^{\sstar}(\E)}
    + \SemiNorm{v}_{\W^{1,\s}(\E)}  \big)
    \qquad\qquad\qquad\forall v \in \W^{1,\s}(\E).
\end{equation}
\end{lemma}
Finally, we recall the existence of a
right-inverse of the trace operator;
cf., e.g., \cite[Chapter~3]{Ern-Guermond:2021}.
\begin{lemma}[Right-inverse of the trace operator]
Let~$\E$ be an element in~$\taun$.
Given $\T$ in~$\TfrakE$
satisfying~\eqref{eq:mesh-regularity-extra}
and~$s$ in~$[1,\infty)$,
there exists a positive constant~$\CtildeTR$
independent of~$\hE$ such that,
for all $v$ in $\W^{\frac{1}{\sprime},\s}(\partial\T)$,
there exists a lifting still denoted by~$v$,
which belongs to~$\W^{1,\s}(\T)$,
such that its trace over~$\partial\T$
coincides with $v$ in~$\W^{\frac{1}{\sprime},\s}(\partial\T)$
and satisfies
\begin{equation} \label{eq:rightinverse-trace}
\Norm{v}_{\W^{1,\s}(\T)}
\le \CtildeTR (\s,\gamma)
    \Norm{v}_{\W^{\frac{1}{\sprime},\s}(\partial\T)} .
\end{equation}
\end{lemma}

The reason why estimate~\eqref{eq:rightinverse-trace}
is stated on simplices
under assumption~\eqref{eq:mesh-regularity-extra}
rather than general elements
satisfying~\eqref{eq:mesh-regularity}
is that in the latter case we are neither aware
nor able to prove bounds with explicit dependence on~$\gamma$.
This will have an impact
in Theorems~\ref{theorem:sharp-ineq}
and~\ref{theorem:nonsharp-ineq} below,
where some estimates only require~\eqref{eq:mesh-regularity},
whereas others need~\eqref{eq:mesh-regularity-extra}.

\subsection{Hessian-based \texorpdfstring{$\W^{2,\p}$}{W²,p} Sobolev--Poincar\'e inequalities} \label{subsection:trivial-inequalities}
Recall the two following Sobolev--Poincar\'e inequalities
from \cite[Theorems~1.7 and 1.8]{Botti-Mascotto:2026}.

\begin{theorem} \label{theorem:BM-W1p}
Let~$\Omega$ be a domain satisfying~ either~\eqref{eq:assumption-Omega-1}
or~\eqref{eq:assumption-Omega-2}.
There exist positive constants $\CBMoo$ and $\CBMot$
such that, for all~$v$ in $\W^{1,\p}(\taun)$,
\small\begin{subequations}
\begin{align}
\Norm{v}_{\L^{\pstar}(\Omega)}
&\le \CBMoo(\p, \gamma, \Omega,\GammaD)
    \Big[ \Norm{\nablaboldh v}_{\Lbf^{\p}(\Omega)}
    + \big(\sum_{\F\in\FcalnID}
        \Norm{\jump{v}}^{\p}_{\L^{\psharp}(\F)}
        \big)^\frac1{\p}  \Big]
        \qquad\text{ if } \p\in[1,d)
        \label{eq:BM-W1p-inequality-1}\\
\Norm{v}_{\L^{\postar}(\Omega)}
&\le \CBMot(\p, \gamma, \Omega,\GammaD)
    \Big[ \Norm{\nablaboldh v}_{\Lbf^{\p}(\Omega)}
    + \big(\sum_{\F\in\FcalnID}
        \htildeF^{1-\p}
        \Norm{\jump{v}}^{\p}_{\L^{\p}(\F)}
        \big)^\frac1{\p}  \Big]
        \qquad\forall \p \in [1,\infty).
        \label{eq:BM-W1p-inequality-2}
\end{align}
\end{subequations}\normalsize
\end{theorem}

Consequences of Theorem~\ref{theorem:BM-W1p}
are the following $\W^{2,\p}$ Sobolev--Poincar\'e inequalities.

\begin{corollary} \label{corollary:trivial-W2p-inequalities}
Let~$\Omega$ be a domain satisfying
either~\eqref{eq:assumption-Omega-1}
or~\eqref{eq:assumption-Omega-2}.
There exist positive constants~$\CBMto$, $\CBMtt$, and~$\CBMtth$
such that, for all~$v$ in $\W^{2,\p}(\taun)$,
\small\begin{subequations}
\begin{align}
\Norm{v}_{\L^{\pstarstar}(\Omega)}
+ \Norm{\nablabold v}_{\Lbf^{\pstar}(\Omega)}
\le & \CBMto(\p, \gamma, \Omega,\GammaD)
   \Big[ \Norm{\Dtwoh v}_{\Lbb^{\p}(\Omega)}
    + \big(\sum_{\F\in\FcalnID}
        \Norm{\jump{\nablaboldh v}}^{\p}_{\Lbf^{\psharp}(\F)}
        \big)^\frac1{\p}\notag\\
&   + \big(\sum_{\F\in\FcalnID}
        \Norm{\jump{v}}^{\pstar}_{\L^{\pstarsharp}(\F)}
        \big)^\frac1{\pstar}
        \Big] \quad\quad \forall \p\in\big[1,\frac d2\big),
        \label{eq:BM-W2p-inequality-1}\\
\Norm{v}_{\L^{\pstarostar}(\Omega)}
+ \Norm{\nablabold v}_{\Lbf^{\pstar}(\Omega)}
\le & \CBMtt(\p, \gamma, \Omega,\GammaD)
   \Big[ \Norm{\Dtwoh v}_{\Lbb^{\p}(\Omega)}
    + \big(\sum_{\F\in\FcalnID}
        \Norm{\jump{\nablaboldh v}}^{\p}_{\Lbf^{\psharp}(\F)}
        \big)^\frac1{\p}\notag\\
&   + \big(\sum_{\F\in\FcalnID}
        \htildeF^{1-\pstar}
        \Norm{\jump{v}}^{\pstar}_{\L^{\pstar}(\F)}
        \big)^\frac1{\pstar}
        \Big] \quad\quad \forall \p\in[1,d),
        \label{eq:BM-W2p-inequality-2}\\
\Norm{v}_{\L^{\postarostar}(\Omega)}
+ \Norm{\nablabold v}_{\Lbf^{\postar}(\Omega)}
\le & \CBMtth(\p, \gamma, \Omega,\GammaD)
 \Big[ \Norm{\Dtwoh v}_{\Lbb^{\p}(\Omega)}
    + \big(\sum_{\F\in\FcalnID} \htildeF^{1-\p}
        \Norm{\jump{\nablaboldh v}}^{\p}_{\Lbf^{\p}(\F)}
        \big)^\frac1{\p} \notag\\
&   + \big(\sum_{\F\in\FcalnID}
        \htildeF^{1-\postar}
        \Norm{\jump{v}}^{\postar}_{\L^{\postar}(\F)}
        \big)^\frac1{\postar}
    \Big] \quad\quad \forall \p\in[1,\infty),
        \label{eq:BM-W2p-inequality-3}\\
\Norm{v}_{\L^{\postar}(\Omega)}
+ \Norm{\nablabold v}_{\Lbf^{\postar}(\Omega)}
\le & \CBMtf(\p, \gamma, \Omega,\GammaD)
 \Big[ \Norm{\Dtwoh v}_{\Lbb^{\p}(\Omega)}
    + \big(\sum_{\F\in\FcalnID} \htildeF^{1-\p}
        \Norm{\jump{\nablaboldh v}}^{\p}_{\Lbf^{\p}(\F)}
        \big)^\frac1{\p} \notag\\
&   + \big(\sum_{\F\in\FcalnID}
        \htildeF^{1-\p}
        \Norm{\jump{v}}^{\p}_{\L^{\p}(\F)}
        \big)^\frac1{\p}
    \Big] \quad\quad \forall \p\in[1,\infty).
        \label{eq:BM-W2p-inequality-4}
\end{align}
\end{subequations}\normalsize
\end{corollary}
\begin{proof}
Using the vector version of~\eqref{eq:BM-W1p-inequality-1}, we write
\small\begin{equation} \label{eq:auxiliary-simple-corollary}
\Norm{\nablabold v}_{\Lbf^{\pstar}(\Omega)}
\le \CBMoo(\p, \gamma, \Omega,\GammaD)
    \Big[ \Norm{\Dtwoh v}_{\Lbb^{\p}(\Omega)}
    + \big(\sum_{\F\in\FcalnID}
        \Norm{\jump{\nablaboldh v}}^{\p}_{\Lbf^{\psharp}(\F)}
        \big)^\frac1{\p}  \Big]
        \qquad\forall \p \in [1,d).
\end{equation}\normalsize
This implies~\eqref{eq:BM-W2p-inequality-1}
since~\eqref{eq:BM-W1p-inequality-1} entails
\[
\Norm{v}_{\L^{\pstarstar}(\Omega)}
\le \CBMoo(\p, \gamma, \Omega,\GammaD)
    \Big[ \Norm{\nablaboldh v}_{\Lbf^{\pstar}(\Omega)}
    + \big(\sum_{\F\in\FcalnID}
        \Norm{\jump{v}}^{\pstar}_{\L^{\pstarsharp}(\F)}
        \big)^\frac1{\pstar}  \Big]
        \qquad\forall\p\in\big[1,\frac d2\big).
\]
Instead, \eqref{eq:auxiliary-simple-corollary}
implies~\eqref{eq:BM-W2p-inequality-2}
since~\eqref{eq:BM-W1p-inequality-2} entails
\[
\Norm{v}_{\L^{\pstarostar}(\Omega)}
\le \CBMot(\p, \gamma, \Omega,\GammaD)
    \Big[ \Norm{\nablaboldh v}_{\Lbf^{\pstar}(\Omega)}
    + \big(\sum_{\F\in\FcalnID}
        \htildeF^{1-\pstar}
        \Norm{\jump{v}}^{\pstar}_{\L^{\pstar}(\F)}
        \big)^\frac1{\pstar}  \Big]
        \qquad\forall \p \in [1,d) .
\]
Using the vector version of~\eqref{eq:BM-W1p-inequality-2}, we write
\small\[
\Norm{\nablabold v}_{\Lbf^{\postar}(\Omega)}
\le \CBMot(\p, \gamma, \Omega,\GammaD)
    \Big[ \Norm{\Dtwoh v}_{\Lbb^{\p}(\Omega)}
    + \big(\sum_{\F\in\FcalnID} \htildeF^{1-\p}
        \Norm{\jump{\nablaboldh v}}^{\p}_{\Lbf^{\p}(\F)}
        \big)^\frac1{\p}  \Big]
        \qquad\forall\p\in[1,\infty).
\]\normalsize
This implies~\eqref{eq:BM-W2p-inequality-3}
since~\eqref{eq:BM-W1p-inequality-2} entails
\small\[
\Norm{v}_{\L^{\postarostar}(\Omega)}
\le \CBMot(\p, \gamma, \Omega,\GammaD)
    \Big[ \Norm{\nablaboldh v}_{\Lbf^{\postar}(\Omega)}
    + \big(\sum_{\F\in\FcalnID}
        \htildeF^{1-\postar}
        \Norm{\jump{v}}^{\postar}_{\L^{\postar}(\F)}
        \big)^\frac1{\postar}  \Big]
        \qquad\forall \p \in [1,\infty).
\]\normalsize
Inequality~\eqref{eq:BM-W2p-inequality-4}
follows analogously.
\end{proof}

\begin{remark} \label{remark:W1p-version-projection}
Assuming that the sequence of meshes has ``no small facets''
and that~$\vh$ is any piecewise polynomial of degree
at most~$\k$ over~$\taun$,
by means of polynomial inverse estimates
and~\eqref{eq:BM-W2p-inequality-4}, we deduce
\[
\begin{split}
\Norm{\vh}_{\L^{\p}(\Omega)}
+ \Norm{\nablaboldh \vh}_{\Lbf^{\p}(\Omega)}
\le & C(\p,\gamma,\Omega,\GammaD)
 \Big[ \Norm{\Dtwoh \vh}_{\Lbb^{\p}(\Omega)}
    + \big(\sum_{\F\in\FcalnID} \hF^{1-\p}
        \Norm{\jump{\nablaboldh \vh}_{N}}^{\p}_{\Lbf^{\p}(\F)}
        \big)^\frac1{\p} \notag\\
&   + \big(\sum_{\F\in\FcalnID}
        \hF^{1-2\p}
        \Norm{\jump{\vh}}^{\p}_{\L^{\p}(\F)}
        \big)^\frac1{\p}
    \Big] \qquad\qquad \forall \p\in[1,\infty).
\end{split}
\]
Furthermore, following \cite[Corollaries 1.10 and 1.11]{Botti-Mascotto:2026},
we also have
\[
\begin{split}
\Norm{\vh}_{\L^{\p}(\Omega)}
+ \Norm{\nablaboldh \vh}_{\Lbf^{\p}(\Omega)}
\le & C(\p,\gamma,\Omega,\GammaD)
 \Big[ \Norm{\Dtwoh \vh}_{\Lbb^{\p}(\Omega)}
    + \big(\sum_{\F\in\FcalnID} \hF^{1-\p}
        \Norm{\PizF \jump{\nablaboldh \vh}_{N}}^{\p}_{\Lbf^{\p}(\F)}
        \big)^\frac1{\p} \notag\\
&   + \big(\sum_{\F\in\FcalnID}
        \hF^{1-2\p}
        \Norm{\PizF \jump{\vh}}^{\p}_{\L^{\p}(\F)}
        \big)^\frac1{\p}
    \Big] \qquad\qquad \forall \p\in[1,\infty).
\end{split}
\]
This inequality is the $\L^{\p}$-version
of the broken Poincar\'e-Friedrichs
in \cite[eq. (1.4)]{Brenner-Wang-Zhao:2004},
in fact with an improved scaling on the standard jump term.
\eremk
\end{remark}

\begin{remark}
The interpretation of inequality~\eqref{eq:BM-W1p-inequality-1}
is rather simple: the indices $\pstar$ and $\psharp$
correspond to those of the continuous embeddings
$\W^{1,\p}(\Omega) \hookrightarrow \L^{\pstar}(\Omega)$
and $\W^{\frac1{\pprime},\p}(\partial\Omega) \hookrightarrow \L^{\psharp}(\partial\Omega)$.
On the other hand, inequality~\eqref{eq:BM-W1p-inequality-2}
involves the index~$\postar$, which has a less immediate
interpretation. We give a heuristic reason
why this index appears.
The crucial question is: what is the maximal, global Sobolev
regularity for the broken space~$\W^{1,\p}(\taun)$?
The answer is $\W^{\s,\p}(\Omega)$ with $\s\p<1$.
In particular, the space $\W^{\frac{1}{\p},\p}(\Omega)$
is a borderline case, where the inclusion
$\W^{1,\p}(\taun) \hookrightarrow \W^{\frac{1}{\p},\p}(\Omega)$
does not hold. The standard Sobolev embedding
$\W^{\frac{1}{\p},\p}(\Omega) \hookrightarrow \L^{\postar}(\Omega)$
heuristically motivates the Lebesgue exponent
in~\eqref{eq:BM-W1p-inequality-2}.
\eremk
\end{remark}

\subsection{Applications of \texorpdfstring{$\W^{1,\p}$}{W1,p} and Hessian-based \texorpdfstring{$\W^{2,\p}$}{W2,p} Sobolev--Poincar\'e inequalities}
\label{subsection:application-1}
In this section, we present concrete applications
of the inequalities established in Theorem~\ref{theorem:BM-W1p}
and Corollary~\ref{corollary:trivial-W2p-inequalities},
focusing on the stability analysis
of nonconforming methods for second-
and fourth-order problems.
In Section~\ref{subsubsection:application-pLaplace},
we consider Crouzeix-Raviart discretizations of
the $\p$-Laplace problem,
while in Section~\ref{subsubsection:application-biharmonic-HH},
we consider nonconforming
discretizations of the biharmonic problem
in Hessian-Hessian formulation.
In what follows, $\k$ denotes a polynomial degree,
$\{\taun\}$ a sequence of either shape-regular
simplicial or regular (in the sense of
Section~\ref{subsection:meshes-broken})
polytopic meshes; $\{\VDG\} = \{\VDG(\taun)\}$
the sequence of spaces of piecewise polynomials
of degree at most~$\k$ over~$\{\taun\}$.

\subsubsection{Nonconforming methods
for the \texorpdfstring{$\p$}{p}-Laplace problem}
\label{subsubsection:application-pLaplace}
Given~$f:\Omega\to\Rbb$ with regularity to be specified later
and $\p\in(1,\infty)$,
consider the $\p$-Laplace problem: find $u:\Omega\to\Rbb$ such that
\[
-\div(\vert\nablabold u\vert^{\p-2} \nablabold u) = f
\text{ in } \Omega,
\qquad\qquad\qquad
u = 0 \text{ on } \Gamma.
\]
Given
\[
\ap(u,v):= (\vert\nablabold u\vert^{\p-2} \nablabold u , 
    \nablabold v)_{\L^2(\Omega)}
\qquad\qquad
\forall u, v \in \W^{1,\p}_0(\Omega),
\]
a weak formulation reads
\begin{equation} \label{eq:pLaplace-weak}
\begin{cases}
\text{find } u \in \W^{1,\p}_0(\Omega) \text{ such that}\\
\ap(u,v) = \langle f, v\rangle
\qquad \forall v \in \W^{1,\p}_0(\Omega).
\end{cases}
\end{equation}
Problem~\eqref{eq:pLaplace-weak}
is well-posed
assuming $f\in [\W^{1,\p}_0(\Omega)]^*=:\W^{-1,\p}(\Omega)$.
Due to the embedding
$\W^{1,\p}(\Omega) \hookrightarrow \L^{\pstar}(\Omega)$,
$f$ in $\L^{\pstarprime}(\Omega)$ guarantees
the well-posedness of the right-hand side
of~\eqref{eq:pLaplace-weak};
we wish to reproduce this condition
on the discrete level with nonconforming methods,
where, unless one uses smoothing operators
on the test function, $f\in \W^{-1,\p}(\Omega)$
is not enough to guarantee the well-posedness
of the discrete right-hand side.

\paragraph*{A Crouzeix-Raviart formulation on triangular meshes.}
We focus here on a Crouzeix-Raviart formulation of problem~\eqref{eq:pLaplace-weak}
on shape-regular triangular meshes $\{\taun\}$.
Introduce the Crouzeix-Raviart space
\[
\VCR:=
\Big\{ \vCR \in \VDG \ \Big|
                \ \int_{\F} \vCR \ \qkmoF = 0 \quad
                \forall \qkmoF \in \Pbb_{\k-1}(\F)
                \quad\forall\F\in\Fcaln \Big\},
\]
the nonlinear form
\[
\apCR(\uCR,\vCR) := (\vert\nablaboldh \uCR \vert^{\p-2} \nablaboldh \uCR,
    \nablaboldh \vCR )_{\L^2(\Omega)}
\qquad\qquad
\forall\uCR,\vCR\in\VCR,
\]
and the positive functional
\[
\Norm{\vCR}_{\apCR}^{\p}
:= \Norm{\nablaboldh \vCR}_{\Lbf^{\p}(\Omega)}^{\p}
    \qquad\qquad \forall\vCR\in\VCR.
\]
Based on the arguments of Remark~\ref{remark:W1p-version-projection},
Theorem~\ref{theorem:BM-W1p},
and the definition of Crouzeix-Raviart spaces,
the functional~$\Norm{\cdot}_{\apCR}$ is a norm on~$\VCR$
and there exists a positive constant $C$ such that
\begin{equation} \label{eq:bound-CR-pstar}
\Norm{\vCR}_{\L^{\pstar}(\Omega)}
    \le C \Norm{\vCR}_{\apCR}
\qquad\qquad
    \forall \vCR\in\VCR,\ \forall \taun \in \{\taun\}.
\end{equation}
Consider the method
\begin{equation} \label{eq:pLaplace-CR}
\begin{cases}
\text{find } \uCR \in \VCR \text{ such that}\\
\apCR(\uCR, \vCR) = (f, \vCR)_{\L^2(\Omega)}
\qquad \forall \vCR \in \VCR.
\end{cases}
\end{equation}
Error estimates for the lowest order
method~\eqref{eq:pLaplace-CR}
were first given in \cite{Liu-Yan:2001} under sufficient
regularity assumptions on the solution to~\eqref{eq:pLaplace-weak};
more recently, error estimates under minimal regularity of the
exact solution were derived in \cite{Kaltenbach:2024, Storn:2026}.

We address the following question:
without introducing further modifications or variational crimes
on the right-hand side in~\eqref{eq:pLaplace-CR},
what regularity on~$f$ guarantees that the method
is well-posed and that stability estimates
with respect to~$f$ can be established?
Assuming $\p\in(1,2)$, the arguments of Remark~\ref{remark:W1p-version-projection}
and Theorem~\ref{theorem:BM-W1p} entail
\small\[
\Norm{\uCR}_{\apCR}^{\p}
= \aCR(\uCR,\uCR)
\overset{\eqref{eq:pLaplace-CR}}{=}
    (f,\uCR)_{\L^2(\Omega)}
\le \Norm{f}_{\L^{\pstarprime}(\Omega)} \Norm{\uCR}_{\L^{\pstar}(\Omega)}
\overset{\eqref{eq:bound-CR-pstar}}{\le}
    C \Norm{f}_{\L^{\pstarprime}(\Omega)} \Norm{\uCR}_{\apCR}.
\]\normalsize
We deduce the estimate
\[
\Norm{\uCR}_{\apCR}
\le C  \Norm{f}_{\L^{\pstarprime}(\Omega)}^{\frac{1}{\p-1}}.
\]
If instead~$\p\in[2,\infty)$,
for~$\q$ arbitrarily large
we have $2\q/(2+\q)<2$ and~$(2\q/(2+\q))^*=\q$, whence
\[
(f,\uCR)_{\L^2(\Omega)}
\le \Norm{f}_{\L^{\qprime}(\Omega)}
        \Norm{\vCR}_{\L^{\q}(\Omega)}
\overset{\eqref{eq:bound-CR-pstar}}{\le}
    C \Norm{f}_{\L^{\qprime}(\Omega)}
    \Norm{\nablaboldh \vCR}_{\Lbf^{\frac{2\q}{2+\q}}(\Omega)}
\le C \Norm{f}_{\L^{\qprime}(\Omega)}
        \Norm{\vCR}_{\apCR}.
\]
Therefore, we get
\[
\Norm{\uCR}_{\apCR}
\le C  \Norm{f}_{\L^{\qprime}(\Omega)}^{\frac{1}{\p-1}}
\qquad\qquad\qquad
\forall\q\in(1,\infty).
\]
Analogous considerations are valid in three dimensions;
the details are omitted.

\subsubsection{Nonconforming methods for the biharmonic problem in Hessian-Hessian formulation}
\label{subsubsection:application-biharmonic-HH}
Given~$f:\Omega\to\Rbb$ with regularity to be specified later,
consider the biharmonic problem: find $u:\Omega\to\Rbb$ such that
\begin{equation} \label{eq:biharmonic-strong:mixed}
\begin{cases}
\Laplacetwo u = f                     &   \text{in } \Omega\\
u=0=\nbfOmega\cdot\nablabold u      &   \text{on } \GammaD\\
\Delta u=0=
\nbfOmega\cdot\nablabold \Delta u   &   \text{on } \GammaN.
\end{cases}
\end{equation}
When~$\GammaN=\emptyset$, given
\[
\btwo(u,v) := (\Dtwo u, \Dtwo v)_{\L^2(\Omega)}
\qquad\qquad\qquad \forall u, v \in \H^2_0(\Omega),
\]
a weak formulation of~\eqref{eq:biharmonic-strong:mixed} reads
\begin{equation} \label{eq:biharmonic-weak:Dirichlet-HH}
\begin{cases}
\text{find } u \in \H^2_0(\Omega) \text{ such that}\\
\btwo(u,v) = \langle f,v \rangle \qquad \forall v\in \H^2_0(\Omega).
\end{cases}
\end{equation}
Problem~\eqref{eq:biharmonic-weak:Dirichlet-HH}
is well-posed thanks to the Lax--Milgram lemma
and standard Poincar\'e inequalities in $\H^2_0(\Omega)$
assuming $f\in [\H^2_0(\Omega)]^*=:\H^{-2}(\Omega)$.

\paragraph*{A symmetric interior penalty Hessian-Hessian formulation.}
First we focus on a symmetric interior penalty discontinuous Galerkin
scheme for problem~\eqref{eq:biharmonic-weak:Dirichlet-HH},
in particular focusing on the Hessian-Hessian formulation;
cf. \cite{Dong-Mascotto-Sutton:2021}.
We consider sequences~$\{\taun\}$ of either shape-regular simplicial
or regular (in the sense of Section~\ref{subsection:meshes-broken})
polytopic meshes; in the latter case, the method
stems combining those in~\cite{Dong:2019,Dong-Mascotto-Sutton:2021}.

Given the space~$\VDG$ of piecewise polynomials of degree
at most~$\k\ge2$ over~$\taun$,
and given the DG penalization parameters
($\Csigma$ and $\Ctau$ are sufficiently large, positive constants
depending on the regularity parameter~$\gamma$
of the sequence of meshes
and the polynomial degree~$\k$)
\begin{equation} \label{eq:biharmonic-DG-parameters}
\sigma := \Csigma \h^{-3},
\qquad\qquad\qquad
\tau   := \Ctau \h^{-1},
\end{equation}
we introduce the bilinear form
\[
\begin{split}
\btwoDG(\uDG,\vDG)
&   := (\Dtwoh \uDG, \Dtwoh \vDG)_{\L^2(\Omega)} \\
&   \quad
    + (\avg{\nbfFcaln \cdot (\divbfh\Dtwoh\uDG)}, \jump{\vDG})_{\Fcaln} 
        - (\avg{(\Dtwoh \uDG)\nbfFcaln}, \jump{\nablaboldh\vDG})_{\Fcaln}\\
&   \quad
    + (\jump{\uDG},\avg{\nbfFcaln \cdot (\divbfh \Dtwoh\vDG)})_{\Fcaln} 
    - (\jump{\nablaboldh\uDG}, \avg{(\Dtwoh \vDG)\nbfFcaln})_{\Fcaln}  \\
&   \quad
    + (\sigma \jump{\uDG}, \jump{\vDG})_{\Fcaln}
    + (\tau \jump{\nablaboldh\uDG}, \jump{\nablaboldh\vDG})_{\Fcaln}
    \qquad \forall \uDG,\vDG\in\VDG
\end{split}
\]
and the positive functional
\[
\Norm{\vDG}_{\btwoDG}^2
:=  \Norm{\Dtwoh \vDG}^2_{\L^2(\Omega)}
    + \Norm{\h^{-\frac12} \jump{\nablaboldh \vDG}}_{\Fcaln}^2
    + \Norm{\h^{-\frac32} \jump{\vDG}}_{\Fcaln}^2
\qquad\qquad\forall\vDG\in\VDG.
\]
Using~\eqref{eq:BM-W2p-inequality-4} entails that
the functional~$\Norm{\cdot}_{\btwoDG}$ is a norm on~$\VDG$
and there exists a positive constant $C$
depending on~$\gamma$ such that
\begin{equation} \label{eq:bound-DG-biharmonic-HH}
\Norm{\vDG}_{\L^2(\Omega)}
\le C \Norm{\vDG}_{\btwoDG}
\qquad\qquad\qquad
    \forall \vDG\in\VDG,\ \forall \taun \in \{\taun\}.
\end{equation}
A symmetric interior penalty discontinuous Galerkin method
for~\eqref{eq:biharmonic-weak:Dirichlet-HH} reads,
cf. \cite[Section~2.3]{Dong-Mascotto-Sutton:2021},
\begin{equation} \label{eq:biharmonic-DG:Dirichlet-HH}
\begin{cases}
\text{find } \uDG \in \VDG \text{ such that}\\
\btwoDG(\uDG, \vDG)
= (f,\vDG)_{\L^2(\Omega)} \qquad \forall \vDG \in \VDG.
\end{cases}
\end{equation}
The following coercivity result holds true,
cf. \cite[Lemma~2.2]{Dong-Mascotto-Sutton:2021}:
there exists a positive constant~$C$
only depending on~$\gamma$
and the polynomial degree~$\k$ such that,
for sufficiently large~$\Csigma$ and~$\Ctau$,
\begin{equation} \label{eq:coercivity-DG-biharmonic-HH}
\Norm{\vDG}_{\btwoDG}^2
\le C \btwoDG(\vDG,\vDG)
\qquad\qquad\qquad
\forall \vDG\in\VDG.
\end{equation}
We address the following question:
without introducing further modifications or variational crimes
on the right-hand side in~\eqref{eq:biharmonic-DG:Dirichlet-HH},
what regularity on~$f$ guarantees that the method
is well-posed and that stability estimates
with respect to~$f$ can be established?
We write, with $\p=2$ in~\eqref{eq:BM-W2p-inequality-4},
\[
\begin{split}
\Norm{\uDG}_{\btwoDG}^2
&   \overset{\eqref{eq:coercivity-DG-biharmonic-HH}}{\le}
        C \btwoDG(\uDG,\uDG)
    \overset{\eqref{eq:biharmonic-DG:Dirichlet-HH}}{=}
        C (f,\uDG)_{\L^2(\Omega)}
    \le C  \Norm{f}_{\L^{\postarprime}(\Omega)}  \Norm{\uDG}_{\L^{\postar}(\Omega)}\\
&   \overset{(\p=2)}{=}
        C \Norm{f}_{\L^{\frac{2 d}{d + 1}}(\Omega)}
        \Norm{\uDG}_{\L^{\frac{2 d}{d - 1}}(\Omega)}
    \overset{\eqref{eq:BM-W2p-inequality-4}}{\le}
        C \CBMtf \Norm{f}_{\L^{\frac{2 d}{d+1}}(\Omega)}
            \Norm{\uDG}_{\btwoDG}.
\end{split}
\]
For a positive constant~$C$ depending
on~$\gamma$ and the polynomial degree~$\k$, we deduce
\[
\Norm{\uDG}_{\btwoDG}
\le C 
\begin{cases}
\Norm{f}_{\L^{\frac43}(\Omega)} &    \text{in 2D,}\\
\Norm{f}_{\L^{\frac32}(\Omega)} &    \text{in 3D}.
\end{cases}
\]
Analogously, with $\p=2$ and $d=3$,
using~\eqref{eq:BM-W2p-inequality-2}
together with polynomial inverse estimates
instead of~\eqref{eq:BM-W2p-inequality-4}
and noting that in this case $\pstarostarprime=9/8$,
we have the improved bound
\[
\Norm{\uDG}_{\btwoDG}
\le C \Norm{f}_{\L^{\frac98}(\Omega)} .
\]
Importantly, the stability bounds above are valid both
in the standard SIPDG case as well as in its polytopic
counterpart as in~\cite{Dong:2019}.

\paragraph*{A $\Ccalz$-interior penalty formulation on simplicial meshes.}
Next, we focus on a $\Ccalz$-interior penalty
(IP) discretization
of~\eqref{eq:biharmonic-weak:Dirichlet-HH}
on simplicial meshes; cf. \cite{Brenner-Sung:2005}.
Compared to the SIPDG formulation in~\eqref{eq:biharmonic-DG:Dirichlet-HH},
instead of~$\VDG$,
we use here the sequence of spaces~$\{\VCcalz\}$
of continuous, piecewise polynomials over~$\taun$
that have zero trace over~$\Gamma$.
For $\tau$ as in~\eqref{eq:biharmonic-DG-parameters},
the $\Ccalz$-IP bilinear form reads
\[
\begin{split}
\btwoCcalz(\uCcalz,\vCcalz)
&   := (\Dtwoh \uCcalz, \Dtwoh \vCcalz)_{\L^2(\Omega)} 
     - (\avg{\nbfFcaln\cdot(\Dtwoh \uCcalz)\nbfFcaln}, \jump{\nablaboldh\vCcalz}_N)_{\Fcaln} \\
&   \quad
    - (\jump{\nablaboldh\uCcalz}_N,
        \avg{\nbfFcaln\cdot(\Dtwoh \vCcalz)\nbfFcaln})_{\Fcaln} \\
&   \quad+ (\tau \jump{\nablaboldh\uCcalz}_N,
        \jump{\nablaboldh\vCcalz}_N)_{\Fcaln}
    \qquad \forall \uCcalz,\vCcalz\in\VCcalz.
\end{split}
\]
We further consider the positive functional
\[
\Norm{\vCcalz}_{\btwoCcalz}^2
:=  \Norm{\Dtwoh \vCcalz}^2_{\L^2(\Omega)}
    + \Norm{\h^{-\frac12}
        \jump{\nablaboldh \vCcalz}_N}_{\Fcaln}^2
\qquad\qquad\forall\vCcalz\in\VCcalz.
\]
The $\Ccalz$-IP version of bounds~\eqref{eq:bound-DG-biharmonic-HH}
and~\eqref{eq:coercivity-DG-biharmonic-HH} is clearly valid
also here:
for sufficiently large~$\Ctau$
and up to using polynomial inverse estimates,
there exists a positive constant $C$
depending on the regularity parameter~$\gamma$
of the sequence of meshes and the polynomial
degree~$\k$ such that
\begin{equation} \label{eq:coercivity-C0IP-biharmonic-HH}
\Norm{\vCcalz}_{\L^2(\Omega)}
\le C \Norm{\vCcalz}_{\btwoCcalz},
\qquad
\Norm{\vCcalz}_{\btwoCcalz}^2
\le C \btwoCcalz(\vCcalz,\vCcalz)
\qquad\qquad\qquad
\forall \vCcalz\in\VCcalz.
\end{equation}
We deduce that~$\Norm{\cdot}_{\btwoCcalz}$
is a norm on~$\VCcalz$.

A $\Ccalz$-IP discontinuous Galerkin method
for~\eqref{eq:biharmonic-weak:Dirichlet-HH} reads,
cf. \cite[Section~4]{Brenner-Sung:2005},
\begin{equation} \label{eq:biharmonic-Ccalz:Dirichlet}
\begin{cases}
\text{find } \uCcalz \in \VCcalz \text{ such that}\\
\btwoCcalz(\uCcalz, \vCcalz)
= (f,\vCcalz)_{\L^2(\Omega)} \qquad \forall \vCcalz \in \VCcalz.
\end{cases}
\end{equation}
We address the following question:
without introducing further modifications or variational crimes
on the right-hand side in~\eqref{eq:biharmonic-Ccalz:Dirichlet},
what regularity on~$f$ guarantees that the method
is well-posed and that stability estimates
with respect to~$f$ can be established?
Proceeding similarly to the case of the SIPDG case~\eqref{eq:biharmonic-DG:Dirichlet-HH},
we write, with $\p=2$ in~\eqref{eq:BM-W2p-inequality-3},
\[
\begin{split}
\Norm{\uCcalz}_{\btwoCcalz}^2
&   \overset{\eqref{eq:coercivity-C0IP-biharmonic-HH}}{\le}
        C \btwoCcalz(\uCcalz,\uCcalz)
    \overset{\eqref{eq:biharmonic-Ccalz:Dirichlet}}{=}
        C (f,\uCcalz)_{\L^2(\Omega)}
    \le C  \Norm{f}_{\L^{\postarostarprime}(\Omega)}  \Norm{\uCcalz}_{\L^{\postarostar}(\Omega)}\\
&   \overset{(\p=2)}{=}
        C \Norm{f}_{\L^{\frac{2 d^2}{d^2 + 2d -1}}(\Omega)}
        \Norm{\uCcalz}_{\L^{\frac{2 d^2}{(d - 1)^2}}(\Omega)}
    \overset{\eqref{eq:BM-W2p-inequality-3}}{\le}
        C \Norm{f}_{\L^{\frac{2 d^2}{d^2 + 2d -1}}(\Omega)}
            \Norm{\uCcalz}_{\btwoCcalz}.
\end{split}
\]
For a positive constant~$C$ depending
on~$\gamma$ and the polynomial degree~$\k$, we deduce
\[
\Norm{\uCcalz}_{\btwoCcalz}
\le C 
\begin{cases}
\Norm{f}_{\L^{\frac87}(\Omega)} &    \text{in 2D},\\
\Norm{f}_{\L^{\frac97}(\Omega)} &    \text{in 3D}.
\end{cases}
\]

\paragraph*{A Morley-Wang-Xu formulation on simplicial meshes.}
Finally, we focus on Morley-Wang-Xu schemes
on shape-regular simplicial meshes
for problem~\eqref{eq:biharmonic-weak:Dirichlet-HH};
cf. \cite{Morley:1968, Wang-Xu:2007}.
Given the Morley-Wang-Xu spaces ($\k=2$)
\begin{equation} \label{eq:MorleyWangXu-definition}
\begin{split}
\VMcal:= \Big\{ \vMcal \in \VDG 
    & \ \Big|\ 
    \int_{\e} (\vMcal{}_{|\E_1}-\vMcal{}_{|\E_2}) d\e=0\;
    \text{ $\forall$ $(d-2)$ dimensional facets~$\e$},\\
    & \quad\quad\quad\quad
        \text{and $\forall$ $\E_1\ne\E_2$
        such that~$\e\subset\partial\E_1\cap\partial\E_2$};\\
    & \quad \int_{\e} \vMcal d\e = 0\;
        \text{ $\forall$ $(d-2)$ dimensional facets~$\e\subset\Gamma$};\\
    & \quad \int_{\F} \jump{\nablaboldh\vMcal}_{N} d\F = 0\;
                \forall \F \in \Fcaln  \Big\},
\end{split}
\end{equation}
we introduce the bilinear form
\[
\btwoMcal(\uMcal,\vMcal) = (\Dtwoh\uMcal,\Dtwoh\vMcal)_{\L^2(\Omega)}
\qquad\qquad
\uMcal,\vMcal \in \VMcal
\]
and the positive functional
\[
\Norm{\vMcal}_{\btwoMcal}^2
:=  \Norm{\Dtwoh \vMcal}^2_{\L^2(\Omega)}
\qquad\qquad\qquad\forall\vMcal\in\VMcal.
\]
In two dimensions, the first condition
in~\eqref{eq:MorleyWangXu-definition} boils down to continuity
at internal vertices and zero value at boundary vertices;
in that case, the space is known as Morley space
while in three dimension it is known as Morley-Wang-Xu space.

A consequence of the definition of~$\VMcal$
is that the functional $\Norm{\cdot}_{\btwoMcal}$ is a norm on~$\VMcal$.
A Morley-Wang-Xu method for~\eqref{eq:biharmonic-weak:Dirichlet-HH} reads,
cf. \cite{Morley:1968},
\begin{equation} \label{eq:biharmonic-Morley:Dirichlet}
\begin{cases}
\text{find } \uMcal \in \VMcal \text{ such that}\\
\btwoMcal(\uMcal, \vMcal)
= (f,\vMcal)_{\L^2(\Omega)} \qquad \forall \vMcal \in \VMcal .
\end{cases}
\end{equation}
We address the following question:
without introducing further modifications or variational crimes
on the right-hand side in~\eqref{eq:biharmonic-Morley:Dirichlet},
what regularity on~$f$ guarantees that the method
is well-posed and that stability estimates
with respect to~$f$ can be established?
We write
\[
\Norm{\uMcal}_{\btwoMcal}^2
= \btwoMcal(\uMcal,\uMcal)
\overset{\eqref{eq:biharmonic-Morley:Dirichlet}}{=}
    (f,\uMcal)_{\L^2(\Omega)}
\le 
    \Norm{f}_{\L^{\frac{2 d^2}{d^2 + 2d - 1}}(\Omega)}
    \Norm{\uMcal}_{\L^{\frac{2 d^2}{(d - 1)^2}}(\Omega)}.
\]
Using~\eqref{eq:BM-W2p-inequality-3},
for a positive constant~$C$ depending
on the regularity parameter~$\gamma$
of the sequence of meshes
and the polynomial degree~$\k$, we deduce
\[
\Norm{\uMcal}_{\btwoMcal}
\le C 
\begin{cases}
\Norm{f}_{\L^{\frac87}(\Omega)} &    \text{in 2D},\\
\Norm{f}_{\L^{\frac97}(\Omega)} &    \text{in 3D}.
\end{cases}
\]

\section{Laplace-based \texorpdfstring{$\W^{2,\p}$}{W²,p} Sobolev--Poincar\'e inequalities} \label{section:main}
In this section, we prove the main results of this work, i.e.,
two types of ``improved'' Sobolev--Poincar\'e inequalities
in $\W^{2,\p}(\taun)$.
The first type of inequality for broken spaces reads as follows.

\begin{theorem} \label{theorem:sharp-ineq}
Assume that~$\Omega$ and~$\p$ are such that the regularity
shift in~\eqref{eq:regularity-shift-Helmholtz} below
is valid for $\q=\pstar$, $\p$ to be fixed below.
Let~$\CHo$, $\CHt$, $\CBcal$, $\CPS$, $\CTR$, and~$\CSob$
be the constants in~\eqref{eq:regularity-shift-Helmholtz} below,
\eqref{eq:bounds:Helmholtz-2} below,
\eqref{eq:stability:RI-curl} below,
\eqref{eq:standard-Poincare-Steklov},
\eqref{eq:local-Sobolev-trace},
and~\eqref{eq:GNS-inequality}, respectively.
Consider~$v$ in $\W^{2,\p}(\taun)$.
Then, further assuming~\eqref{eq:assumption-Omega-3}
on~$\Omega$, we have
\begin{equation} \label{eq:sharp-W2p-inequalities-1}
\begin{split}
\Norm{\nablaboldh v}_{\Lbf^{\pstar}(\Omega)} 
&    \le \CtildeBMoo(\p,\gamma,\Omega,\GammaD)
        \Norm{\Deltah v}_{\L^{\p}(\Omega)}\\
&   \quad   + \CtildeBMot(\p,\gamma,\Omega,\GammaD)
      \big( \sum_{\F \in \FcalnID} 
        \Norm{\jump{\nablaboldh v}_{N}}_{\L^{\psharp}(\F)}^{\p} \big)^\frac1{\p}\\
&   \quad + \CtildeBMoth(\p,\gamma,\Omega,\GammaD)
        \big( \sum_{\F \in \FcalnIN} 
        \Norm{\jump{\nablaboldh v}_{T}}_{\Lbf^{\psharp}(\F)}^{\p}   \big)^\frac1{\p}
        \qquad\forall\p\in[1,d),
\end{split}
\end{equation}\normalsize
where
\small\[
\begin{split}
&   \CtildeBMoo(\p,\gamma,\Omega,\GammaD)
    := \CHo(\pstarprime,\Omega,\GammaD)
    \CSob(\pstarprime,\Omega)
    (1+\CPS(\pstarprime)),\\
&   \CtildeBMot(\p,\gamma,\Omega,\GammaD)
:= \text{C}_{\text{TR}}^{\frac{\pstar}{\pprime}}
    (\pstarprime,\gamma)
    2^{\frac{1}{\p}}
        \big( 1+\CSob(\pstarprime,\Omega)
        (1+\CPS(\pstarprime)) \big) 
        \CHo(\pstarprime,\Omega,\GammaD),\\
&   \CtildeBMoth(\p,\gamma,\Omega,\GammaD)
:= \text{C}_{\text{TR}}^{\frac{\pstar}{\pprime}} 
    (\pstarprime,\gamma)
    2^{\frac{1}{\p}}\cdot\\
&   \qquad\qquad\qquad\qquad\qquad
        \cdot \big( 1+\CSob(\pstarprime,\Omega)
        (1+\CPS(\pstarprime)) \big) 
        \CHt(\pstarprime,\Omega,\GammaD)
        \CBcal(\pstarprime,\Omega,\GammaD).
\end{split}
\]\normalsize
Moreover, given~$\CtildeTR$ and~$\CEG$
as in~\eqref{eq:rightinverse-trace}
and~\eqref{eq:stability-zero-extension} below,
for domains~$\Omega$ with general topology
and further assuming~\eqref{eq:mesh-regularity-extra}
on the sequence of meshes,
\small\begin{equation} \label{eq:sharp-W2p-inequalities-2}
\begin{split}
\Norm{\nablaboldh v}_{\Lbf^{\pstar}(\Omega)}
&   \le  \CtildeBMto(\p,\gamma,\Omega,\GammaD)
        \Norm{\Deltah v}_{\L^{\p}(\Omega)}
     + \CtildeBMtt(\p,\gamma,\Omega,\GammaD)
        \big( \sum_{\F \in \FcalnID} 
        \Norm{\jump{\nablaboldh v}_{N}}_{\L^{\psharp}(\F)}^{\p} \big)^\frac1{\p}\\
&   \quad + \CtildeBMtth(\p,\gamma,\Omega,\GammaD)
    \big( \sum_{\F \in \FcalnIN} {}
    \Norm{\jump{v} }_{\W^{\frac{1}{\pstarprime},\pstar}(\F)}^{\pstar}
    \big)^\frac{1}{\pstar}
    \qquad \forall\p\in
    \Big[1, \frac{2d}{d+2} \Big) ,
\end{split}
\end{equation}\normalsize
where
\[
\begin{split}
&   \CtildeBMto(\p,\gamma,\Omega,\GammaD):= \CtildeBMoo(\p,\gamma,\Omega,\GammaD),
    \quad
    \CtildeBMtt(\p,\gamma,\Omega,\GammaD):= \CtildeBMot(\p,\gamma,\Omega,\GammaD),\\
&   \CtildeBMtth(\p,\gamma,\Omega,\GammaD)
    := \CtildeTR(\pstar,\gamma)
        \CHt(\pstarprime,\Omega,\GammaD)
        \CEG(\pstar,\gamma).
\end{split}
\]
\end{theorem}

The other inequality is a variant
with nonoptimal Lebesgue indices
but which allows for more general $\p$.
\begin{theorem} \label{theorem:nonsharp-ineq}
Assume that~$\Omega$ and~$\p$ are such that the regularity
shift in~\eqref{eq:regularity-shift-Helmholtz} below
is valid for $\q=\pstar$, $\p$ to be fixed below.
Consider~$v$ in $\W^{2,\p}(\taun)$.
Let~$\CHo$, $\CHt$, $\CBcal$, $\CPS$,
$\CSob$, and~$\CMT$
be the constants in~\eqref{eq:regularity-shift-Helmholtz} below,
\eqref{eq:bounds:Helmholtz-2} below,
\eqref{eq:stability:RI-curl} below,
\eqref{eq:standard-Poincare-Steklov},
\eqref{eq:GNS-inequality},
and~\eqref{eq:CMT-definition} below, respectively.
Then, further assuming~\eqref{eq:assumption-Omega-3}
on $\Omega$, we have
\begin{equation} \label{eq:nonsharp-W2p-inequalities-1}
\begin{split}
\Norm{\nablaboldh v}_{\Lbf^{\postar}(\Omega)} 
&    \le \CtildeBMtho(\p,\gamma,\Omega,\GammaD)
        \Norm{\Deltah v}_{\L^{\p}(\Omega)}\\
&   \quad + \CtildeBMtht(\p,\gamma,\Omega,\GammaD)
      \Big( \sum_{\F \in \FcalnID}{}
            \htildeF^{1-\p} \Norm{\jump{\nablaboldh v}_N}_{\L^{\p}(\F)}^{\p} \Big)^\frac{1}{\p}\\
&   \quad + \CtildeBMthth(\p,\gamma,\Omega,\GammaD)
        \Big( \sum_{\F \in \FcalnIN}{}
            \htildeF^{1-\p} \Norm{\jump{\nablaboldh v}_T}_{\Lbf^{\p}(\F)}^{\p} \Big)^\frac{1}{\p}
            \quad\forall \p\in[1,\infty),
\end{split}
\end{equation}
where
\[
\begin{split}
&   \CtildeBMtho(\p,\gamma,\Omega,\GammaD)
    := \CHo(\postarprime,\Omega,\GammaD)
        \vert\Omega\vert^{\frac{d-1}{d}}
                \CSob(\postarprime,\Omega)
                (1+\CPS(\postarprime)),\\
&   \CtildeBMtht(\p,\gamma,\Omega,\GammaD)
    := \CHo(\postarprime,\Omega,\GammaD)
        \CMT(\p,\Omega), \\
&   \CtildeBMthth(\p,\gamma,\Omega,\GammaD)
    := \CMT(\p,\Omega)
        \CBcal(\postarprime,\Omega,\GammaD)
        \CHt(\postarprime,\Omega,\GammaD).
\end{split}
\]
Moreover,
given~$\CtildeTR$ and~$\CEG$
as in~\eqref{eq:rightinverse-trace}
and~\eqref{eq:stability-zero-extension} below,
for domains~$\Omega$ with general topology
and further assuming~\eqref{eq:mesh-regularity-extra}
on the sequence of meshes,
\small\begin{equation} \label{eq:nonsharp-W2p-inequalities-2}
\begin{split}
\Norm{\nablaboldh v}_{\Lbf^{\postar}(\Omega)} 
&    \le \CtildeBMfo(\p,\gamma,\Omega,\GammaD)
        \Norm{\Deltah v}_{\L^{\p}(\Omega)}\\
&   \quad    + \CtildeBMft(\p,\gamma,\Omega,\GammaD)
      \Big( \sum_{\F \in \FcalnID}{}
            \htildeF^{1-\p} \Norm{\jump{\nablaboldh v}_N}_{\L^{\p}(\F)}^{\p} \Big)^\frac{1}{\p}\\
&   \quad + \CtildeBMfth(\p,\gamma,\Omega,\GammaD)
        \big( \sum_{\F \in \FcalnIN} {}
    \Norm{\jump{v} }_{\W^{\frac{1}{\postarprime},\postar}(\F)}^{\postar}
    \big)^\frac{1}{\postar}
    \quad\forall \p\in \big[1,\frac{2(d-1)}{d}\big),
\end{split}
\end{equation}\normalsize
where
\[
\begin{split}
&\CtildeBMfo(\p,\gamma,\Omega,\GammaD):= \CtildeBMtho(\p,\gamma,\Omega,\GammaD),
\qquad
\CtildeBMft(\p,\gamma,\Omega,\GammaD):= \CtildeBMtht(\p,\gamma,\Omega,\GammaD),\\
&
\CtildeBMfth(\p,\gamma,\Omega,\GammaD)
            := \CMT(\p,\Omega)
            \CBcal(\postarprime,\Omega,\GammaD)
            \CHt(\postarprime,\Omega,\GammaD).
\end{split}
\]
\end{theorem}

While applications of these inequalities
are exhibited in Section~\ref{subsection:application-2},
the proofs of Theorems~\ref{theorem:sharp-ineq}
and~\ref{theorem:nonsharp-ineq}
are given in Sections~\ref{subsection:sharp-inequalities}
and~\ref{subsection:nonsharp-inequalities}, respectively.
The proofs are detailed in three dimensions;
the corresponding proof in two dimensions is analogous.
A crucial tool in the analysis is the Helmholtz
decomposition, which is discussed in Section~\ref{subsection:Helmholtz-dec}.

\subsection{Application of Laplace-based \texorpdfstring{$\W^{2,\p}$}{W²,p} Sobolev--Poincar\'e inequalities} \label{subsection:application-2}
In this section, we present applications
of the inequalities established
in Theorems~\ref{theorem:sharp-ineq}
and~\ref{theorem:nonsharp-ineq},
focusing on the stability analysis
of nonconforming methods for fourth-order problems.
In Sections~\ref{subsubsection:application-pBiharmonic}
and~\ref{subsubsection:application-biharmonic-LL},
we consider nonconforming discretizations
of the $\p$-biharmonic problem
and of the biharmonic problem
in Laplace-Laplace formulation
with mixed boundary conditions, respectively.
In what follows, $\k$ denotes a polynomial degree,
$\{\taun\}$ a sequence of either shape-regular
simplicial or regular (in the sense of
Section~\ref{subsection:meshes-broken})
polytopic meshes;
$\VDG$ the corresponding sequence of spaces
of piecewise polynomials of degree at most~$\k$.

\subsubsection{Nonconforming methods for
the \texorpdfstring{$\p$}{p}-biharmonic problem} \label{subsubsection:application-pBiharmonic}
Given~$\p\in(1,\infty)$ and $f:\Omega\to\Rbb$ with regularity to be specified later,
consider the $\p$-biharmonic problem: find $u:\Omega\to\Rbb$ such that
\begin{equation} \label{eq:pbiharmonic-strong:mixed}
\begin{cases}
\Delta(\vert\Delta u\vert^{\p-2} \Delta u) = f
        &   \text{in } \Omega\\
u=0=\nbfOmega\cdot \nablabold u
        &   \text{on } \Gamma.
\end{cases}
\end{equation}
Given
\[
\cp(u,v) := (\vert \Delta u \vert^{\p-2} \Delta u, \Delta v)_{\L^2(\Omega)}
\qquad\qquad\qquad \forall u, v \in \W^{2,\p}_0(\Omega),
\]
a weak formulation of~\eqref{eq:pbiharmonic-strong:mixed} reads
\begin{equation} \label{eq:pbiharmonic-weak:Dirichlet}
\begin{cases}
\text{find } u \in \W^{2,\p}_0(\Omega) \text{ such that}\\
\cp(u,v) = \langle f,v\rangle \qquad \forall v\in \W^{2,\p}_0(\Omega).
\end{cases}
\end{equation}
Problem~\eqref{eq:pbiharmonic-weak:Dirichlet}
is well-posed
assuming $f\in [\W^{2,\p}_0(\Omega)]^*=\W^{-2,\p}(\Omega)$.
If $\p=2$, then formulation~\eqref{eq:pbiharmonic-weak:Dirichlet}
is equivalent to the ``Hessian-Hessian'' formulation
\[
\begin{cases}
\text{find } u \in \W^{2,\p}_0(\Omega) \text{ such that}\\
(\vert \Dtwo u \vert^{\p-2} \Dtwo u, \Dtwo v)_{\L^2(\Omega)}
= \langle f,v\rangle \qquad \forall v\in \W^{2,\p}_0(\Omega).
\end{cases}
\]
A $\Ccalz$-IP discretization of this problem is studied in~\cite{Hoppe:2022}.

\paragraph*{A Morley-Wang-Xu formulation on triangular meshes.}
We focus here on a Morley-Wang-Xu
scheme for problem~\eqref{eq:pbiharmonic-weak:Dirichlet}
on shape-regular simplicial meshes;
cf. \cite{Morley:1968,Wang-Xu:2007}.
Throughout, we assume that~$\p$ must satisfy the regularity
shift in~\eqref{eq:regularity-shift-Helmholtz}
and also $\p \in (1,d)$;
see Remark~\ref{remark:application-p-bounds}
below for details.

We consider sequences~$\{\taun\}$ of shape-regular simplicial meshes.
Given the Morley-Wang-Xu space~$\VMcal$
in~\eqref{eq:MorleyWangXu-definition},
introduce the bilinear form
\[
\begin{split}
\cpMcal(\uMcal,\vMcal)
&   = (\vert \Deltah\uMcal \vert^{\p-2} \Deltah\uMcal,
        \Deltah\vMcal)_{\L^2(\Omega)}\\
&   \quad 
    + \sum_{\F\in\Fcaln} \htildeF^{1-\p}
        (\vert \jump{\nablaboldh \uMcal}_{N}\vert^{\p-2}
        \jump{\nablaboldh \uMcal}_{N},
        \jump{\nablaboldh \vMcal}_{N})_{\F}
        \\
&   \quad
    + \sum_{\F\in\Fcaln} \htildeF^{1-\p}
        (\vert\jump{\nablaboldh\uMcal}_{T}\vert^{\p-2}
        \jump{\nablaboldh \uMcal}_{T},
        \jump{\nablaboldh \vMcal}_{T})_{\L^2(\F)}
    \qquad
    \uMcal,\vMcal \in \VMcal
\end{split}
\]
and the positive functional
\small\[
\begin{split}
\Norm{\vMcal}_{\cpMcal}^{\p}
&   :=  \Norm{\Deltah \vMcal}_{\L^{\p}(\Omega)}^{\p}
    + \sum_{\F\in\Fcaln} \htildeF^{1-\p}
    \Norm{\jump{\nablaboldh \vMcal}_{N}}_{\L^{\p}(\F)}^{\p}
    + \sum_{\F\in\Fcaln} \htildeF^{1-\p}
    \Norm{\jump{\nablaboldh\vMcal}_T}_{\Lbf^{\p}(\F)}^{\p}
\quad\forall\vMcal\in\VMcal.
\end{split}
\]\normalsize
The functional~$\Norm{\cdot}_{\cpMcal}$
is a norm on~$\VMcal$.
In fact, \eqref{eq:nonsharp-W2p-inequalities-1} implies
the existence of a positive constant~$C$
only depending on the regularity parameter~$\gamma$
of the sequence of meshes such that
\[
\Norm{\nablaboldh \vMcal}_{\Lbf^{\p}(\Omega)}
\le C \Norm{\vMcal}_{\cpMcal}
\qquad\qquad
\forall\p\in(1,d).
\]
On the other hand,
the functional on the left-hand side above
is a norm on~$\VMcal$: when it is zero,
$\vMcal$ is piecewise constant over~$\taun$;
exploiting the continuity of~$\vMcal$
on the $(d-2)$ dimensional facets of~$\taun$
and the fact that~$\vMcal$ is zero
on the $(d-2)$ dimensional boundary facets,
this implies that~$\vMcal=0$.
A Morley-Wang-Xu method ($\k=2$)
for~\eqref{eq:pbiharmonic-weak:Dirichlet} reads
\begin{equation} \label{eq:pbiharmonic-Morley:Dirichlet}
\begin{cases}
\text{find } \uMcal \in \VMcal \text{ such that}\\
\cpMcal(\uMcal, \vMcal)
= (f,\vMcal)_{\L^2(\Omega)} \qquad \forall \vMcal \in \VMcal .
\end{cases}
\end{equation}
We address the following question:
without introducing further modifications or variational crimes
on the right-hand side in~\eqref{eq:pbiharmonic-Morley:Dirichlet},
what regularity on~$f$ guarantees that the method
is well-posed and that stability estimates
with respect to~$f$ can be established?
Observe that also the regularity
shift in~\eqref{eq:regularity-shift-Helmholtz}
must be taken into account.
We write
\[
\Norm{\uMcal}_{\cpMcal}^{\p}
= \cpMcal(\uMcal,\uMcal)
\overset{\eqref{eq:pbiharmonic-Morley:Dirichlet}}{=}
    (f,\uMcal)_{\L^2(\Omega)}
\le 
    \Norm{f}_{\L^{\pstarostarprime}(\Omega)}
    \Norm{\uMcal}_{\L^{\pstarostar}(\Omega)}.
\]
We focus on the last term on the right-hand side
and 
then, for a positive constant~$C$, we have
\[
\begin{split}
&   \Norm{\uMcal}_{\L^{\pstarostar}(\Omega)}
   \overset{\eqref{eq:BM-W1p-inequality-1}}{\le}
    C \big[ \Norm{\nablaboldh \uMcal}_{\Lbf^{\pstar}(\Omega)}
        + \big(\sum_{\F\in\Fcaln} \htildeF^{1-\pstar}
        \Norm{\jump{\uMcal}}_{\L^{\pstar}(\F)}
        ^{\pstar} \big)^{\frac{1}{\pstar}}  \big]\\
&   \overset{\eqref{eq:sharp-W2p-inequalities-1}}{\le}
    C \Big[ \Norm{\Deltah \uMcal}_{\L^{\p}(\Omega)}
        + \big(\sum_{\F\in\Fcaln}
        \Norm{\jump{\nablaboldh \uMcal}_N}_{\L^{\psharp}(\F)}
        ^{\p} \big)^{\frac{1}{\p}}\\
&   \qquad\qquad
        + \big(\sum_{\F\in\FcalnI}
        \Norm{\jump{\nablaboldh\uMcal}_T}_{\Lbf^{\psharp}(\F)}
        ^{\p} \big)^{\frac{1}{\p}}
        + \big(\sum_{\F\in\Fcaln} \htildeF^{1-\pstar}
        \Norm{\jump{\uMcal}}_{\L^{\pstar}(\F)}
        ^{\pstar} \big)^{\frac{1}{\pstar}}
        \Big].
\end{split}
\]
Using a polynomial inverse estimate,
a Poincar\'e-type estimates
with the properties of Morley-Wang-Xu space,
and~\eqref{eq:mesh-regularity-extra},
we write
\[
\begin{split}
\htildeF^{1-\pstar} \Norm{\jump{\uMcal}}_{\L^{\pstar}(\F)}^{\pstar}
&   \le C \htildeF^{1-\pstar +(d-1)\big(1-\frac{\pstar}{\psharp}\big)}
        \Norm{\jump{\uMcal}}_{\L^{\psharp}(\F)}^{\pstar}\\
&   \le C \htildeF^{1 +(d-1)\big(1-\frac{\pstar}{\psharp}\big)}
    \Norm{\jump{\nablaboldh\uMcal}_T}_{\Lbf^{\psharp}(\F)}^{\pstar}
    \qquad\forall\F\in\Fcaln.
\end{split}
\]
We deduce
\[
\begin{split}
\big(\sum_{\F\in\Fcaln} \htildeF^{1-\pstar}
        \Norm{\jump{\uMcal}}_{\L^{\pstar}(\F)}
        ^{\pstar} \big)^{\frac{1}{\pstar}}
&   \le C
    \Big(\sum_{\F\in\Fcaln}
    \htildeF^{1 +(d-1)\big(1-\frac{\pstar}{\psharp}\big)}
       \Norm{\jump{\nablaboldh\uMcal}_T}_{\Lbf^{\psharp}(\F)}^{\pstar} \Big)^{\frac{\p}{\pstar}\frac{1}{\p}}\\
&   \le \Big(\sum_{\F\in\Fcaln}
    \htildeF^{\big[1
    +(d-1)\big(1-\frac{\pstar}{\psharp}\big)\big]\frac{\p}{\pstar}}
       \Norm{\jump{\nablaboldh\uMcal}_T}_{\Lbf^{\psharp}(\F)}^{\p} \Big)^{\frac{1}{\p}}.
\end{split}
\]
A computation reveals
\[
\big[1 +(d-1)\big(1-\frac{\pstar}{\psharp}\big)\big]
    \frac{\p}{\pstar} = 0.
\]
This, together with polynomial inverse estimates, implies
\[
\begin{split}
\Norm{\uMcal}_{\L^{\pstarostar}(\Omega)}
&   \le C \Big[ \Norm{\Deltah \uMcal}_{\L^{\p}(\Omega)}
        + \big(\sum_{\F\in\Fcaln}
        \Norm{\jump{\nablaboldh \uMcal}_N}_{\L^{\psharp}(\F)}
        ^{\p} \big)^{\frac{1}{\p}}
        + \big(\sum_{\F\in\Fcaln}
        \Norm{\jump{\nablaboldh\uMcal}_T}_{\Lbf^{\psharp}(\F)}
        ^{\p} \big)^{\frac{1}{\p}}\Big]\\
&   \le C \Big[ \Norm{\Deltah \uMcal}_{\L^{\p}(\Omega)}
        + \big(\sum_{\F\in\Fcaln} \htildeF^{1-\p}
        \Norm{\jump{\nablaboldh \uMcal}_N}_{\L^{\p}(\F)}
        ^{\p} \big)^{\frac{1}{\p}}\\
&   \qquad\qquad
        + \big(\sum_{\F\in\Fcaln} \htildeF^{1-\p}
        \Norm{\jump{\nablaboldh\uMcal}_T}_{\Lbf^{\p}(\F)}
        ^{\p} \big)^{\frac{1}{\p}}\Big]
        =: C \Norm{\uMcal}_{\cpMcal}.
\end{split}
\]
We deduce
\begin{equation} \label{eq:regime-2}
\Norm{\uMcal}_{\cpMcal}
\le C
    \Norm{f}_{\L^{\pstarostarprime}(\Omega)}^{\frac{1}{\p-1}}
    \qquad\qquad\qquad \forall \p\in (1,d)
    \text{ satisfying~\eqref{eq:regularity-shift-Helmholtz} below}.
\end{equation}

\begin{remark} \label{remark:application-p-bounds}
For the application of Theorems~\ref{theorem:sharp-ineq}
and~\ref{theorem:nonsharp-ineq},
we have to guarantee that the regularity
shift in~\eqref{eq:regularity-shift-Helmholtz}
below holds true.
Based on the arguments
in Appendix~\ref{appendix:regularity-Poisson},
if full Dirichlet boundary conditions
are considered, cf.~\eqref{eq:JK-bounds},
one has that the admissible range of~$\p$ is
\begin{equation} \label{eq:Savare}
\p \in \Big(\frac{2d}{d+1},\frac{2d}{d-1}\Big).
\end{equation}
As such, the stability estimate~\eqref{eq:regime-2}
combined with~\eqref{eq:regularity-shift-Helmholtz}
and~\eqref{eq:Savare} reads
\[
\Norm{\uMcal}_{\cpMcal}
\le C \Norm{f}_{\L^{\pstarostarprime}(\Omega)}^{\frac{1}{\p-1}}
\qquad\qquad
\begin{cases}
\forall \p\in \Big( \frac{4}{3}, 2 \Big)
            \text{ if } d=2,\\
\forall \p\in \Big( \frac{3}{2}, 3 \Big)
            \text{ if } d=3.
\end{cases}
\]
\eremk
\end{remark}

\subsubsection{Nonconforming methods for the biharmonic problem in Laplace-Laplace formulation} 
\label{subsubsection:application-biharmonic-LL}
Given~$f:\Omega\to\Rbb$ with regularity to be specified later,
consider the biharmonic problem
\begin{equation} \label{eq:biharmonic-strong:mixed-variant}
\begin{cases}
\Laplacetwo u = f                   &   \text{in } \Omega\\
u=0=\Delta u                        &   \text{on } \GammaN\\
\nbfOmega\cdot\nablabold u =0=
\nbfOmega\cdot\nablabold \Delta u   &   \text{on } \GammaD.
\end{cases}
\end{equation}
We assume that $\GammaN\ne\emptyset\ne\GammaD$;
observe that the subscripts $D$ and $N$ are not related
to the standard notation corresponding
to Dirichlet and Neumann boundaries. Given
\[
\btildetwo(u,v) := (\Delta u, \Delta v)
\qquad \forall u, v \in 
\VLL:= \{ v\in \H^1_{\GammaN}(\Omega) \cap \H^2(\Omega) \mid
    \nbfOmega\cdot\nablabold v=0
    \text{ on }\GammaD\},
\]
a weak formulation of~\eqref{eq:biharmonic-strong:mixed-variant} reads
\begin{equation} \label{eq:biharmonic-weak:Dirichlet-LL}
\begin{cases}
\text{find } u \in \VLL \text{ such that}\\
\btildetwo(u,v) = \langle f,v\rangle
    \qquad \forall v\in \VLL.
\end{cases}
\end{equation}
Problem~\eqref{eq:biharmonic-weak:Dirichlet-LL}
is well-posed thanks to the Lax--Milgram lemma
and standard Poincar\'e inequalities
in $\VLL$
assuming $f\in [\VLL]^*$.

\paragraph*{A symmetric interior penalty Laplace-Laplace formulation.}
We focus on a symmetric interior penalty discontinuous Galerkin
scheme in Laplace-Laplace formulation
for problem~\eqref{eq:biharmonic-weak:Dirichlet-LL};
cf., e.g., \cite{Mozolevski-Suli-Bosing:2007, Georgoulis-Houston:2009}.
We consider sequences~$\{\taun\}$ of either shape-regular simplicial
or regular (in the sense of Section~\ref{subsection:meshes-broken})
polytopic meshes; in the latter case, the method
is that in~\cite{Dong:2019}.

Given the space~$\VDG$ of piecewise polynomials of degree
at most~$\k\ge2$ over~$\taun$,
and given DG penalization parameters
($\Csigma$ and $\Ctau$ are sufficiently large, positive constants
depending on the regularity parameter~$\gamma$
of the sequence of meshes and the polynomial degree~$\k$)
as in~\eqref{eq:biharmonic-DG-parameters},
introduce the bilinear form
\small\[
\begin{split}
\btildetwoDG(\uDG,\vDG)
&   := (\Deltah \uDG, \Deltah \vDG)_{\L^2(\Omega)} \\
&   \quad
    + (\avg{\nbfFcaln \cdot (\nablaboldh\Deltah\uDG)}, \jump{\vDG})_{\FcalnIN} 
        - (\avg{\Deltah \uDG}, \jump{\nablaboldh\vDG}_{N})_{\FcalnID}\\
&   \quad
    + (\jump{\uDG},\avg{\nbfFcaln \cdot (\nablaboldh \Deltah\vDG)})_{\FcalnIN} 
    - (\jump{\nablaboldh\uDG}_{N}, \avg{\Deltah \vDG})_{\FcalnID}  \\
&   \quad
    + (\sigma \jump{\uDG}, \jump{\vDG})_{\FcalnIN}
    + (\tau \jump{\nablaboldh\uDG}_{N},
        \jump{\nablaboldh\vDG}_{N})_{\FcalnID}
    \qquad \forall \uDG,\vDG\in\VDG
\end{split}
\]\normalsize
and the positive functional
\[
\Norm{\vDG}_{\btildetwoDG}^2
:=  \Norm{\Deltah \vDG}^2_{\L^2(\Omega)}
    + \Norm{\h^{-\frac12} \jump{\nablaboldh \vDG}_{N}}_{\FcalnID}^2
    + \Norm{\h^{-\frac32} \jump{\vDG}}_{\FcalnIN}^2
    \qquad\forall\vDG\in\VDG.
\]
A consequence of Theorem~\ref{theorem:BM-W1p}
and Theorem~\ref{theorem:nonsharp-ineq}
is that the functional~$\Norm{\cdot}_{\btildetwoDG}$
is a norm on~$\VDG$ and there exists a positive constant $C$
such that
\[
\Norm{\vDG}_{\L^2(\Omega)}
\le C \Norm{\vDG}_{\btildetwoDG}
\qquad\qquad\qquad
    \forall \vDG \in \{\VDG\},\ \forall \taun \in \{\taun\}.
\]
A symmetric interior penalty discontinuous Galerkin method
for~\eqref{eq:biharmonic-weak:Dirichlet-LL} reads,
cf. \cite[Section~3]{Mozolevski-Suli-Bosing:2007}
and \cite[Section~3]{Georgoulis-Houston:2009},
\begin{equation} \label{eq:biharmonic-DG:Dirichlet-LL}
\begin{cases}
\text{find } \uDG \in \VDG \text{ such that}\\
\btildetwoDG(\uDG, \vDG)
= (f,\vDG)_{\L^2(\Omega)} \qquad \forall \vDG \in \VDG.
\end{cases}
\end{equation}
The following coercivity result holds true,
cf. \cite[Section~5.1]{Georgoulis-Houston:2009}:
there exists a positive constant~$C$
only depending on~$\gamma$
and the polynomial degree~$\k$
such that, for sufficiently large~$\Csigma$ and~$\Ctau$,
\begin{equation} \label{eq:coercivity-DG-biharmonic-LL}
\Norm{\vDG}_{\btildetwoDG}^2
\le C \btildetwoDG(\vDG,\vDG)
\qquad\qquad\qquad
\forall \vDG\in\VDG.
\end{equation}
We address the following question:
without introducing further modifications or variational crimes
on the right-hand side in~\eqref{eq:biharmonic-DG:Dirichlet-LL},
what regularity on~$f$ guarantees that the method
is well-posed and that stability estimates
with respect to~$f$ can be established?
Observe that the regularity
shift in~\eqref{eq:regularity-shift-Helmholtz}
holds true for~$\p=2$;
cf. Remark~\ref{remark:regularity-shift-qprime=2} below.
We write, with $\p=2$ in~\eqref{eq:BM-W1p-inequality-2}
and~\eqref{eq:nonsharp-W2p-inequalities-1},
\[
\begin{split}
\Norm{\uDG}_{\btildetwoDG}^2
&   \overset{\eqref{eq:coercivity-DG-biharmonic-LL}}{\le}
        C \btildetwoDG(\uDG,\uDG)
    \overset{\eqref{eq:biharmonic-DG:Dirichlet-LL}}{=}
        (f,\uDG)_{\L^2(\Omega)}
    \le \Norm{f}_{\L^{\postarostarprime}(\Omega)}  \Norm{\uDG}_{\L^{\postarostar}(\Omega)}\\
 &  \overset{(\p=2)}{=}
        \Norm{f}_{\L^{\frac{2 d^2}{d^2 + 2d - 1}}(\Omega)}
        \Norm{\uDG}_{\L^{\frac{2 d^2}{(d - 1)^2}}(\Omega)}
    \overset{\eqref{eq:BM-W1p-inequality-2},
            \eqref{eq:nonsharp-W2p-inequalities-1}}{\le}
        C \Norm{f}_{\L^{\frac{2 d^2}{d^2 + 2d - 1}}(\Omega)}
            \Norm{\uDG}_{\btildetwoDG},
\end{split}
\]
where in the last step we further use polynomial inverse
estimates on facets and considerations analogous to those
for the Morley-Wang-Xu scheme
in Section~\ref{subsubsection:application-pBiharmonic}.

For a positive constant~$C$ depending
on the regularity parameter~$\gamma$ of the sequence of meshes
and the polynomial degree~$\k$, we deduce
\[
\Norm{\uDG}_{\btildetwoDG}
\le C 
\begin{cases}
\Norm{f}_{\L^{\frac87}(\Omega)} &    \text{in 2D}\\
\Norm{f}_{\L^{\frac97}(\Omega)} &    \text{in 3D}.
\end{cases}
\]
These stability bounds are valid both in the standard SIPDG case
as well as in its polytopic counterpart as in~\cite{Dong:2019}.
Analogous bounds may be derived for other nonconforming
methods including $\Ccalz$-IP methods and Morley-Wang-Xu,
as we did in Section~\ref{subsubsection:application-biharmonic-HH}
for the Hessian-Hessian formulation;
we omit the details.

\subsection{A Helmholtz decomposition} \label{subsection:Helmholtz-dec}
Let~$\q$ be in $(1,\infty)$.
Given $\sigmabold$ in $\Lbf^{\qprime}(\Omega)$,
we have the following Helmholtz decomposition.
There exists $w$ in $\W^{1,\qprime}_{\GammaN}(\Omega)$
and $\xibold$ in $\Lbf^{\qprime}(\Omega)$
satisfying $\div \xibold =0$
and~$\xibold\cdot\nbfOmega=0$ on~$\GammaD$
such that
\begin{equation} \label{eq:Helmholtz-decomposition}
\sigmabold = \nablabold w + \xibold.
\end{equation}
We further introduce the space
\[
\W_{00}^{\frac1{\qprime},\q}(\GammaD)
:= \{ w \in \W^{\frac1{\qprime},\q}(\GammaD)
    \mid w/\rho^{\frac1{\qprime}} \in \L^{\q}(\GammaD) \},
\qquad
\rho(x)=\text{dist}(x,\partial\GammaD)
\quad x \in \GammaD.
\]
We construct explicitly the two terms on the right-hand side
of~\eqref{eq:Helmholtz-decomposition}.
We pick $w$ as the solution to
\begin{equation} \label{eq:1st-element-Helmholtz}
\begin{cases}
-\Delta w = -\div \sigmabold
    & \text{in } [\W^{1,\q}_{\GammaN}(\Omega)]^* \\
\nablabold w \cdot \nbfOmega = \sigmabold \cdot \nbfOmega 
    & \text{in } [\W_{00}^{\frac1{\qprime},\q}(\GammaD)]^*\\
w = 0  
    & \text{on } \GammaN.
\end{cases}
\end{equation}
We assume the well-posedness of problem~\eqref{eq:1st-element-Helmholtz}
and the following a priori estimate:
there exists a positive constant
depending on $\qprime$ and~$\Omega$ such that
\begin{equation} \label{eq:regularity-shift-Helmholtz}
\SemiNorm{w}_{\W^{1,\qprime}(\Omega)}
\le \CHo(\qprime,\Omega,\GammaD)
    \Norm{\sigmabold}_{\Lbf^{\qprime}(\Omega)}.
\end{equation}
\begin{remark} \label{remark:regularity-shift-qprime=2}
The regularity result in~\eqref{eq:regularity-shift-Helmholtz}
is trivial in the Hilbert case~$\qprime=2$
but is not generally true if $\qprime\ne2$.
On the one hand, the right-hand side cannot be
``too singular'', which gives a lower bound on~$\qprime$;
on the other hand, one should take into account
the shape of the domain~$\Omega$, which introduces
singularities in the solution,
yielding an upper bound on~$\qprime$;
cf. the discussion
in Appendix~\ref{appendix:regularity-Poisson}.
\eremk
\end{remark}

Define
\[
\xibold := \sigmabold-\nablabold w \in \Lbf^{\qprime}(\Omega).
\]
Using~\eqref{eq:regularity-shift-Helmholtz},
there exists a positive constant~$\CHt$
depending on $\qprime$ and~$\Omega$ such that
\begin{equation} \label{eq:bounds:Helmholtz-2}
\Norm{\xibold}_{\Lbf^{\qprime}(\Omega)}
    \le \CHt(\qprime,\Omega,\GammaD)
        \Norm{\sigmabold}_{\Lbf^{\qprime}(\Omega)}.
\end{equation}
Note that~$\xibold$ is divergence free
and~$\xibold\cdot\nbfOmega=0$ on~$\GammaD$.
If~$\Omega$ satisfies \eqref{eq:assumption-Omega-3},
then \cite[Theorem~3]{Hiptmair-Pechstein:2020}
and~\cite{Galdi:2011} imply in 3D the existence
of~$\Psibold$ in~$\Wbf^{1,\qprime}_{\GammaD}(\Omega)$
and of a positive constant~$\CBcal$
only depending on~$\qprime$, $\GammaD$, and~$\Omega$
such that
\small\begin{equation} \label{eq:stability:RI-curl}
\curlbold\Psibold=\xibold,
\qquad
    \SemiNorm{\Psibold}_{\W^{1,\qprime}(\Omega)}
    \le \CBcal(\qprime,\Omega,\GammaD)
        \Norm{\xibold}_{\Lbf^{\qprime}(\Omega)}
    \overset{\eqref{eq:bounds:Helmholtz-2}}{\le}
        \CBcal(\qprime,\Omega,\GammaD)  
        \CHt(\qprime,\Omega,\GammaD)
        \Norm{\sigmabold}_{\Lbf^{\qprime}(\Omega)}.
\end{equation}\normalsize
In the 2D case, the regularized right-inverse
of the $\curlbold$ operator
is a scalar function.
If~\eqref{eq:assumption-Omega-1} and $\GammaN=\emptyset$
further hold true, then we also have an explicit upper bound
for~$\CBcal$ in terms of the diameter
and largest inradius of~$\Omega$;
cf. \cite[Corollary 19 (ii)]{Guzman-Salgado:2021}.

\subsection{Proof of the sharp inequalities} \label{subsection:sharp-inequalities}
Here, we prove Theorem~\ref{theorem:sharp-ineq}.

\paragraph*{A duality rewriting.}
Let~$v$ be in~$\W^{2,\p}(\taun)$, $\p \in [1,d)$.
By duality, we have
\begin{equation} \label{eq:initial-step-W2-broken}
\Norm{\nablaboldh v}_{\Lbf^{\pstar}(\Omega)}
= \sup_{\sigmabold \in \Lbf^{\pstarprime}(\Omega)}
    \frac{_{\Lbf^{\pstarprime}} ( \sigmabold, \nablaboldh v )_{\Lbf^{\pstar}}}
        {\Norm{\sigmabold}_{\Lbf^{\pstarprime}(\Omega)}}.
\end{equation}
Using decomposition~\eqref{eq:Helmholtz-decomposition} with~$\q=\pstar$
and assuming~\eqref{eq:regularity-shift-Helmholtz},
we have the splitting
\[
\sigmabold = \nablabold w + \xibold,
\qquad\qquad
w \in \W^{1,\pstarprime}_{\GammaN}(\Omega),
\quad
\xibold \in \Lbf^{\pstarprime}(\Omega).
\]
We rewrite the numerator on the right-hand side
of~\eqref{eq:initial-step-W2-broken}
\begin{equation} \label{eq:initial-split-Helmholtz}
_{\Lbf^{\pstarprime}} ( \sigmabold, \nablaboldh v )_{\Lbf^{\pstar}}
    = _{\Lbf^{\pstarprime}} ( \nablabold w, \nablaboldh v )_{\Lbf^{\pstar}}
    + _{\Lbf^{\pstarprime}} ( \xibold, \nablaboldh v )_{\Lbf^{\pstar}}
    =: T_1 + T_2.
\end{equation}
We handle the two terms on the right-hand side separately.

\paragraph*{The term~$T_1$.}
Integrating by parts, we have
\[
T_1
= _{\L^{\pprime}(\Omega)} (w, \Deltah v )_{\L^{\p}(\Omega)}
    + \sum_{\F \in \FcalnID} {}
        _{\L^{\psharpprime}(\F)} ( w , 
        \jump{\nablaboldh v}_N)_{\L^{\psharp}(\F)}
=:  T_{1,1} + T_{1,2}.
\]
As for the term~$\T_{1,1}$, we use H\"older's inequality,
the identity $\pprime=\pstarprimestar$,
the Sobolev embedding $\W^{1,\pstarprime}(\Omega) \hookrightarrow \L^{\pstarprimestar}(\Omega)$
with a Poincar\'e inequality, and write
\[
\begin{split}
T_{1,1}
&   \le  \Norm{w}_{\L^{\pprime}(\Omega)} \Norm{\Deltah v}_{\L^{\p}(\Omega)}
    = \Norm{w}_{\L^{\pstarprimestar}(\Omega)} \Norm{\Deltah v}_{\L^{\p}(\Omega)}\\
&   \overset{\eqref{eq:standard-Poincare-Steklov}, \eqref{eq:GNS-inequality}}{\le}
    \CSob(\pstarprime,\Omega)(1+\CPS(\pstarprime))
    \SemiNorm{w}_{\W^{1,\pstarprime}(\Omega)} \Norm{\Deltah v}_{\L^{\p}(\Omega)}\\
&   \overset{\eqref{eq:regularity-shift-Helmholtz}}{\le}
    \CSob(\pstarprime,\Omega)(1+\CPS(\pstarprime))
    \CHo(\pstarprime,\Omega,\GammaD)
    \Norm{\sigmabold}_{\Lbf^{\pstarprime}(\Omega)} \Norm{\Deltah v}_{\L^{\p}(\Omega)}.
\end{split}
\]
As for the term~$T_{1,2}$, noting $\pstarprimesharp = \psharpprime$,
$r^*=r^\sharp\ostar$ for all $r$ in $[1,d)$,
and again $\pstarprimestar=\pprime$, we write
\small\begin{equation} \label{eq:crazy-estimates-faces}
\begin{split}
&   \sum_{\F\in\FcalE} \Norm{w}_{\L^{\psharpprime}(\F)}^{\pstarprimestar}
    = \sum_{\F\in\FcalE} \Norm{w}_{\L^{\pstarprimesharp}(\F)}^{\pstarprimestar}
    \overset{\text{(Jensen)}}{\le}
    \big( \sum_{\F\in\FcalE} \Norm{w}_{\L^{\pstarprimesharp}(\F)}^{\pstarprimesharp} \big)^{\ostar}\\
&   = \Norm{w}_{\L^{\pstarprimesharp}(\partial\E)}^{\pstarprimesharp\ostar}
    = \Norm{w}_{\L^{\pstarprimesharp}(\partial\E)}^{\pprime}
    \overset{\eqref{eq:local-Sobolev-trace}}{\le}
    \text{C}_{\text{TR}}^{\pstar}(\pstarprime,\gamma)
    2^{\pprime-1}
    \big(\Norm{w}_{\L^{\pstarprimestar}(\E)}^{\pstarprimestar}
    +\SemiNorm{w}_{\W^{1,\pstarprime}(\E)} ^{\pstarprimestar}\big).
\end{split}
\end{equation}\normalsize
H\"older's inequality for sequences with indices~$\p$ and~$\pprime$ yields
\[
T_{1,2}
    \le \Big( \sum_{\F \in \FcalnID}  
            \Norm{w}_{\L^{\psharpprime}(\F)}^{\pprime} \Big)^\frac{1}{\pprime}
        \Big( \sum_{\F \in \FcalnID} 
            \Norm{\jump{\nablaboldh v}_{N}}_{\L^{\psharp}(\F)}^{\p}   \Big)^\frac1{\p}.
\]
We are left with estimating
the first term on the right-hand side:
\small\[
\begin{split}
&   \Big( \sum_{\F \in \FcalnID}  
        \Norm{w}_{\L^{\psharpprime}(\F)}^{\pprime} \Big)^\frac{1}{\pprime}
    \le   \Big( \sum_{\E \in \taun} \sum_{\F\in\FcalE}
        \Norm{w}_{\L^{\psharpprime}(\F)}^{\pprime} \Big)^\frac{1}{\pprime}\\
&   \overset{\eqref{eq:crazy-estimates-faces}}{\le}
    \text{C}_{\text{TR}}^{\frac{\pstar}{\pprime}}
    (\pstarprime,\gamma)
    2^{\frac{1}{\p}}
        \Big( \sum_{\E \in \taun} \big(
        \Norm{w}_{\L^{\pprime}(\E)}^{\pprime}
            +\SemiNorm{w}_{\W^{1,\pstarprime}(\E)} ^{\pprime} \big)
            \Big)^\frac{1}{\pprime}\\
&   = \text{C}_{\text{TR}}^{\frac{\pstar}{\pprime}}
    (\pstarprime,\gamma)
    2^{\frac{1}{\p}}
    \Big( \Norm{w}_{\L^{\pprime}(\Omega)}^{\pprime}
            +\sum_{\E \in \taun} \SemiNorm{w}_{\W^{1,\pstarprime}(\E)} ^{\pprime}
                \Big)^\frac{1}{\pprime} \\
&   \overset{\text{(Jensen)}}{\le}
    \text{C}_{\text{TR}}^{\frac{\pstar}{\pprime}} 
    (\pstarprime,\gamma)
    2^{\frac{1}{\p}}
    \Big( \Norm{w}_{\L^{\pprime}(\Omega)}^{\pprime}
            + \big(\sum_{\E \in \taun} \SemiNorm{w}_{\W^{1,\pstarprime}(\E)}^{\pstarprime}
            \big)  ^{\frac{\pprime}{\pstarprime}}
                \Big)^\frac{1}{\pprime}\\
&   = \text{C}_{\text{TR}}^{\frac{\pstar}{\pprime}}
    (\pstarprime,\gamma)
    2^{\frac{1}{\p}}
        \Big( \Norm{w}_{\L^{\pprime}(\Omega)}^{\pprime}
            + \SemiNorm{w}_{\W^{1,\pstarprime}(\Omega)}^{\pprime}
                \Big)^\frac{1}{\pprime}
    \le \text{C}_{\text{TR}}^{\frac{\pstar}{\pprime}}
    (\pstarprime,\gamma)
    2^{\frac{1}{\p}}
        \big(\Norm{w}_{\L^{\pprime}(\Omega)}
            + \SemiNorm{w}_{\W^{1,\pstarprime}(\Omega)}\big)\\
&   \overset{\eqref{eq:standard-Poincare-Steklov},\eqref{eq:GNS-inequality}}{\le}
        \text{C}_{\text{TR}}^{\frac{\pstar}{\pprime}}
        (\pstarprime,\gamma)
        2^{\frac{1}{\p}}
        \big( 1+\CSob(\pstarprime,\Omega)
        (1+\CPS(\pstarprime)) \big)
        \SemiNorm{w}_{\W^{1,\pstarprime}(\Omega)}\\
&   \overset{\eqref{eq:regularity-shift-Helmholtz}}{\le}
    \text{C}_{\text{TR}}^{\frac{\pstar}{\pprime}}
    (\pstarprime,\gamma) 2^{\frac{1}{\p}}
        \big( 1+\CSob(\pstarprime,\Omega)
        (1+\CPS(\pstarprime)) \big) 
        \CHo(\pstarprime,\Omega,\GammaD)
        \Norm{\sigmabold}_{\Lbf^{\pstarprime}(\Omega)}.
\end{split}
\]\normalsize
Combining the two displays above we arrive at
\[
\begin{split}
T_{1,2}
&   \le \text{C}_{\text{TR}}^{\frac{\pstar}{\pprime}} 
    (\pstarprime,\gamma)
    2^{\frac{1}{\p}}
        \big( 1+\CSob(\pstarprime,\Omega)
        (1+\CPS(\pstarprime)) \big)   
        \CHo(\pstarprime,\Omega,\GammaD) \cdot \\
&   \quad \cdot \Big( \sum_{\F \in \FcalnID} 
            \Norm{\jump{\nablaboldh v}_{N}}_{\L^{\psharp}(\F)}^{\p}   \Big)^\frac1{\p}
        \Norm{\sigmabold}_{\Lbf^{\pstarprime}(\Omega)}.
\end{split}
\]
Combining the bounds on $T_{1,1}$ and~$T_{1,2}$ yields
\small\begin{equation} \label{eq:T1-version-1&2}
\begin{split}
T_1
&   \le
    \CHo(\pstarprime,\Omega,\GammaD) 
    \Big[ \CSob(\pstarprime,\Omega)
    (1+\CPS(\pstarprime))
        \Norm{\Deltah v}_{\L^{\p}(\Omega)}\\
&   \quad
    +\text{C}_{\text{TR}}^{\frac{\pstar}{\pprime}}
    (\pstarprime,\gamma)
    2^{\frac{1}{\p}}
    \big( 1+\CSob(\pstarprime,\Omega)
    (1+\CPS(\pstarprime)) \big)\cdot\\
&   \qquad\qquad\qquad \cdot  \big( \sum_{\F \in \FcalnID} 
        \Norm{\jump{\nablaboldh v}_{N}}_{\L^{\psharp}(\F)}^{\p} \big)^\frac1{\p} \Big]
        \Norm{\sigmabold}_{\Lbf^{\pstarprime}(\Omega)}.
\end{split}
\end{equation}\normalsize

\paragraph*{The term~$T_2$.}
Owing to~\eqref{eq:assumption-Omega-3},
we can consider the regularized right-inverse~$\Psibold$
in $\Wbf^{1,\pstarprime}(\Omega)$ of the $\curlbold$ operator of~$\xibold$
and write
\[
\begin{split}
T_2
&   = _{\Lbf^{\pstarprime}} (\xibold, \nablaboldh v )_{\Lbf^{\pstar}}
    =  \sum_{\E\in\taun}{}
        _{\Lbf^{\pstarprime}(\E)} (\curlbold \Psibold, \nablaboldh v )_{\Lbf^{\pstar}(\E)}\\
&   \overset{(\pstarprimesharp=\pstarprime)}{=}
    \sum_{\E\in\taun}{}_{\Lbf^{\psharpprime}(\partial\E)}
        (\Psibold, \nbfE \times \nablaboldh v)_{\Lbf^{\psharp}(\partial\E)} 
    = \sum_{\F \in \FcalnIN} {}
        _{\Lbf^{\psharpprime}(\F)} (\Psibold  ,
            \jump{\nablaboldh v}_T)_{\Lbf^{\psharp}(\F)} \\
&   = \sum_{\F \in \FcalnIN} {}
        _{\Lbf^{\pstarprimesharp}(\F)} (\Psibold  ,
            \jump{\nablaboldh v}_T)_{\Lbf^{\pstarprimesharpprime}(\F)}.
\end{split}
\]
Proceeding as for the term~$T_{1,2}$, we write
\begin{equation} \label{eq:crazy-estimates-faces-il-ritorno}
\begin{split}
&   \sum_{\F \in \FcalnIN}
        \Norm{\Psibold}_{\Lbf^{\psharpprime}(\F)}^{\pstarprimestar}
    = \sum_{\F \in \FcalnIN}
        \Norm{\Psibold}_{\Lbf^{\pstarprimesharp}(\F)}^{\pstarprimestar}
    \overset{\text{(Jensen)}}{\le}
        \big( \sum_{\F \in \FcalnIN}
            \Norm{\Psibold}_{\Lbf^{\pstarprimesharp}(\F)}^{\pstarprimesharp} \big)^{\ostar}\\
&   = \big( \Norm{\Psibold}_{\Lbf^{\pstarprimesharp}(\partial\E)} 
            ^{\pstarprimesharp} \big)^{\ostar}
    \overset{\eqref{eq:local-Sobolev-trace}}{\le}
    \text{C}_{\text{TR}}^{\pstar}
    (\pstarprime,\gamma)
    2^{\pprime-1}
    \big(\Norm{\Psibold}_{\Lbf^{\pstarprimestar}(\E)}^{\pstarprimestar}
    +\SemiNorm{\Psibold}_{\Wbf^{1,\pstarprime}(\E)} ^{\pstarprimestar}\big).
\end{split}
\end{equation}
H\"older's inequality for sequences with indices~$\p$ and~$\pprime$ yields
\[
T_2
    \le \Big( \sum_{\F \in \FcalnIN}  
            \Norm{\Psibold}_{\Lbf^{\psharpprime}(\F)}^{\pprime} \Big)^\frac{1}{\pprime}
        \Big( \sum_{\F \in \FcalnIN} 
            \Norm{\jump{\nablaboldh v}_{T}}_{\Lbf^{\psharp}(\F)}^{\p}   \Big)^\frac1{\p}.
\]
Since the second term on the right-hand side will appear in the final estimate,
we only focus on the first one:
\footnotesize\[
\begin{split}
&   \Big( \sum_{\F \in \FcalnIN}
        \Norm{\Psibold}_{\Lbf^{\psharpprime}(\F)}^{\pprime} \Big)^\frac{1}{\pprime}
    \le   \Big( \sum_{\E \in \taun} \sum_{\F\in\FcalE}
        \Norm{\Psibold}_{\Lbf^{\psharpprime}(\F)}^{\pprime} \Big)^\frac{1}{\pprime}\\
&   \overset{\eqref{eq:crazy-estimates-faces-il-ritorno}}{\le}
    \text{C}_{\text{TR}}^{\frac{\pstar}{\pprime}}
    (\pstarprime,\gamma)
    2^{\frac{1}{\p}}
        \Big( \sum_{\E \in \taun} \big(
        \Norm{\Psibold}_{\Lbf^{\pprime}(\E)}^{\pprime}
            +\SemiNorm{\Psibold}_{\Wbf^{1,\pstarprime}(\E)} ^{\pprime} \big)
            \Big)^\frac{1}{\pprime}\\
&   = \text{C}_{\text{TR}}^{\frac{\pstar}{\pprime}}
    (\pstarprime,\gamma)
    2^{\frac{1}{\p}}
    \Big( \Norm{\Psibold}_{\Lbf^{\pprime}(\Omega)}^{\pprime}
            +\sum_{\E \in \taun} \SemiNorm{\Psibold}_{\Wbf^{1,\pstarprime}(\E)} ^{\pprime}
                \Big)^\frac{1}{\pprime} \\
&   \overset{\text{(Jensen)}}{\le}
    \text{C}_{\text{TR}}^{\frac{\pstar}{\pprime}}
    (\pstarprime,\gamma)
    2^{\frac{1}{\p}}
    \Big( \Norm{\Psibold}_{\Lbf^{\pprime}(\Omega)}^{\pprime}
            + \big(\sum_{\E \in \taun} \SemiNorm{\Psibold}_{\Wbf^{1,\pstarprime}(\E)}^{\pstarprime}
            \big)  ^{\frac{\pprime}{\pstarprime}}
                \Big)^\frac{1}{\pprime}\\
&   = \text{C}_{\text{TR}}^{\frac{\pstar}{\pprime}}
    (\pstarprime,\gamma)
    2^{\frac{1}{\p}}
    \Big( \Norm{\Psibold}_{\Lbf^{\pprime}(\Omega)}^{\pprime}
            + \SemiNorm{\Psibold}_{\Wbf^{1,\pstarprime}(\Omega)}^{\pprime}
                \Big)^\frac{1}{\pprime}
    \le
    \text{C}_{\text{TR}}^{\frac{\pstar}{\pprime}}
    (\pstarprime,\gamma)
    2^{\frac{1}{\p}}
    \big(\Norm{\Psibold}_{\Lbf^{\pprime}(\Omega)}
            + \SemiNorm{\Psibold}_{\Wbf^{1,\pstarprime}(\Omega)}\big)\\
&   \overset{\eqref{eq:standard-Poincare-Steklov},\eqref{eq:GNS-inequality}}{\le}
     \text{C}_{\text{TR}}^{\frac{\pstar}{\pprime}}
     (\pstarprime,\gamma)
     2^{\frac{1}{\p}}
        \big( 1+\CSob(\pstarprime,\Omega)
        (1+\CPS(\pstarprime)) \big)
        \SemiNorm{\Psibold}_{\Wbf^{1,\pstarprime}(\Omega)}\\
&   \overset{\eqref{eq:stability:RI-curl}}{\le}
        \text{C}_{\text{TR}}^{\frac{\pstar}{\pprime}}
        (\pstarprime,\gamma)
        2^{\frac{1}{\p}}
        \big( 1+\CSob(\pstarprime,\Omega)
        (1+\CPS(\pstarprime)) \big)
        \CHt(\pstarprime,\Omega,\GammaD)
        \CBcal(\pstarprime,\Omega,\GammaD)
        \Norm{\sigmabold}_{\Lbf^{\pstarprime}(\Omega)}.
\end{split}
\]\normalsize
Combining the two displays above yields
\small\begin{equation}\label{eq:T2-version-1}
\begin{split}
T_2
&   \le \text{C}_{\text{TR}}^{\frac{\pstar}{\pprime}}
    (\pstarprime,\gamma)
    2^{\frac{1}{\p}}
        \big( 1+\CSob(\pstarprime,\Omega)
        (1+\CPS(\pstarprime)) \big)
        \CHt(\pstarprime,\Omega,\GammaD)
        \CBcal(\pstarprime,\Omega,\GammaD) \cdot\\
&   \quad \cdot \big( \sum_{\F \in \FcalnIN} 
        \Norm{\jump{\nablaboldh v}_{T}}_{\Lbf^{\psharp}(\F)}^{\p}   \big)^\frac1{\p}
        \Norm{\sigmabold}_{\Lbf^{\pstarprime}(\Omega)}.
\end{split}
\end{equation}\normalsize

\paragraph*{Wrap up.}
Collecting~\eqref{eq:initial-split-Helmholtz},
\eqref{eq:T1-version-1&2}, and~\eqref{eq:T2-version-1}, we arrive at
\footnotesize\[
\begin{split}
&   \Norm{\nablaboldh v}_{\Lbf^{\pstar}(\Omega)}
\le \CHo(\pstarprime,\Omega,\GammaD) \CSob(\pstarprime,\Omega)
    (1+\CPS(\pstarprime)) \Norm{\Deltah v}_{\L^{\p}(\Omega)}\\
&   \quad
    + \text{C}_{\text{TR}}^{\frac{\pstar}{\pprime}}
    (\pstarprime,\gamma)
    2^{\frac{1}{\p}}
        \big( 1+\CSob(\pstarprime,\Omega)
        (1+\CPS(\pstarprime)) \big) 
        \CHo(\pstarprime,\Omega,\GammaD)
    \big( \sum_{\F \in \FcalnIN} 
        \Norm{\jump{\nablaboldh v}_{T}}_{\Lbf^{\psharp}(\F)}^{\p}   \big)^\frac1{\p}\\
&   \quad +
    \text{C}_{\text{TR}}^{\frac{\pstar}{\pprime}}
    (\pstarprime,\gamma)
    2^{\frac{1}{\p}}
        \big( 1+\CSob(\pstarprime,\Omega)
        (1+\CPS(\pstarprime)) \big) \cdot\\
&   \qquad\qquad\qquad
    \cdot\CHt(\pstarprime,\Omega,\GammaD)
        \CBcal(\pstarprime,\Omega,\GammaD)
    \big( \sum_{\F \in \FcalnID} 
        \Norm{\jump{\nablaboldh v}_{N}}_{\L^{\psharp}(\F)}^{\p} \big)^\frac1{\p},
\end{split}
\]\normalsize
which is~\eqref{eq:sharp-W2p-inequalities-1}.

\paragraph*{The term~$T_2$: alternative approach.}
We now assume that~$\Omega$ has general topology
and $\pstar\in(1,2)$, i.e., $\p\in(1,2d/(d+2))$.
Set for simplicity~$\q:=\pstar$.
For all facets~$\F$ in~$\Fcaln$, let~$\EF$
be the element of~$\taun$
such that~$\F$ belongs to $\mathcal F^{\EF}$
and $\nbfEF=\nbfF$;
moreover, let~$\TEF$ be any simplex in~$\TfrakEF$,
cf. the notation in Section~\ref{subsection:meshes-broken},
such that~$\partial\TEF\cap\partial\EF
=: \FTEF \subseteq \F$.
All such simplices~$\T$ are collected in~$\tautildeEF$.

We denote the zero extension operator from~$\FTEF$
to~$\partial\TEF$ by~$\EfrakTEF$.
The following stability result follows from
\cite[Lemma~17.1]{Ern-Guermond:2021} and a scaling argument
(it is here crucial the fact that we consider the simplex
$\TEF$, otherwise the scaling argument would be valid
with an unknown constant):
there exists a positive constant~$\CEG$
only depending on~$\q$ and~$\gamma$ such that
\begin{equation} \label{eq:stability-zero-extension}
\Norm{\EfrakTEF {\jump{v}}}_{\W^{\frac{1}{\qprime},\q}(\partial\TEF)}
\le \CEG(\q,\gamma)
\Norm{\jump{v}}_{\W^{\frac{1}{\qprime},\q}(\FTEF)} .
\end{equation}
The left-hand side of~\eqref{eq:stability-zero-extension}
is well-posed since $\q<2$ implies $\qprime>2$;
cf. \cite[Chapter 17]{Ern-Guermond:2021}.

Since~$\xibold$ in~\eqref{eq:Helmholtz-decomposition}
is divergence free, an integration by parts implies
\[
\begin{split}
T_2
&   := _{\Lbf^{\qprime}} ( \xibold, \nablaboldh v )_{\Lbf^{\q}}
    = \sum_{\E\in\taun}{}
        _{\Lbf^{\qprime}(\E)} ( \xibold, \nablaboldh v )_{\Lbf^{\q}(\E)}
    = \sum_{\E \in \taun} {}
        _{\W^{-\frac{1}{\qprime},\qprime}(\partial\E)}
        \langle \xibold\cdot\nbfE, 
            v \rangle_{\W^{\frac{1}{\qprime},\q}(\partial\E)}\\
&   = \sum_{\E \in \taun} 
    \sum_{\F\in\FcalE}
    \sum_{\TEF \in \tautildeEF}{}
        _{\W^{-\frac{1}{\qprime},\qprime}(\partial\TEF)}
        \langle \xibold\cdot\nbfE, 
        \EfrakTEF (v_{|\E}) \rangle_{\W^{\frac{1}{\qprime},\q}(\partial\TEF)}.
\end{split}
\]
We deduce
\[
\begin{split}
T_2
&   =    \sum_{\F \in \FcalnIN} 
            \sum_{\TEF \in \tautildeEF}{}
        _{\W^{-\frac{1}{\qprime},\qprime}(\partial\TEF)}
        \langle \xibold\cdot\nbfF , 
        \EfrakTEF\jump{v} \rangle _{\W^{\frac{1}{\qprime},\q}(\partial\TEF)}\\
&   \le \sum_{\F \in \FcalnIN} 
    \sum_{\TEF \in \tautildeEF}{}
    \Norm{\xibold\cdot\nbfF }
    _{\W^{-\frac{1}{\qprime},\qprime}
    (\partial\TEF)}
    \Norm{\EfrakTEF\jump{v}}_{\W^{\frac{1}{\qprime},\q}
    (\partial\TEF)}.
\end{split}
\]
We focus on the first local terms separately:
since~$\xibold$ in~\eqref{eq:Helmholtz-decomposition}
is divergence free and~\eqref{eq:mesh-regularity-extra}
holds true
(it is here crucial again the fact that we consider the simplex
$\TEF$, otherwise \eqref{eq:rightinverse-trace} would be valid
with an unknown constant),
\[
\begin{split}
&   \Norm{\xibold\cdot\nbfF}_{\W^{-\frac{1}{\qprime}, \qprime}(\partial\TEF)}
    :=\sup_{\psi \in \W^{\frac{1}{\qprime},\q}
    (\partial\TEF)}
    \frac{_{\W^{-\frac{1}{\qprime},\qprime}
    (\partial\TEF)}
    \langle \xibold\cdot\nbfF, \psi
    \rangle_{\W^{\frac{1}{\qprime},\q}
    (\partial\TEF)}}
    {\Norm{\psi}_{\W^{\frac{1}{\qprime},\q}
    (\partial\TEF)}} \\
&   \overset{\eqref{eq:rightinverse-trace}}{\le}
    \CtildeTR(\q,\gamma)
    \sup_{\psi\in\W^{1,\q}(\TEF)} 
    \frac{_{\W^{-\frac{1}{\qprime},\qprime}
    (\partial\TEF)}
    \langle \xibold\cdot\nbfF, \psi
    \rangle_{\W^{\frac{1}{\qprime},\q}
    (\partial\TEF)}}
    {\SemiNorm{\psi}_{\W^{1,\q} (\TEF)}}\\
&   = \CtildeTR(\q,\gamma)
    \sup_{\psi\in\W^{1,\q}(\TEF)} 
    \frac{_{\Lbf^{\qprime}(\TEF)}(\xibold,\nabla\psi)_{\Lbf^{\q}(\TEF)}}
    {\SemiNorm{\psi}_{\W^{1,\q} (\TEF)}}
    \le \CtildeTR(\q,\gamma)
    \Norm{\xibold}_{\Lbf^{\qprime}(\TEF)}.
\end{split}
\]
Combining the two displays above,
and using H\"older's inequality for sequences
with indices~$\qprime$ and~$\q$ yield
\begin{equation} \label{eq:T2-version-12}
\begin{split}
T_2
&   \le \CtildeTR(\q,\gamma)
    \Norm{\xibold}_{\Lbf^{\qprime}(\Omega)}
    \Big( \sum_{\F \in \FcalnIN}
    \sum_{\TEF \in \tautildeEF}{}
    \Norm{\EfrakEF\jump{v}}_{\W^{\frac{1}{\qprime},\q}
    (\partial\TEF)}^{\q}
    \Big)^\frac{1}{\q}\\
&   \overset{\eqref{eq:bounds:Helmholtz-2}}{\le}
    \CtildeTR(\q,\gamma)
    \CHt(\qprime,\Omega,\GammaD)
    \Norm{\sigmabold}_{\Lbf^{\qprime}(\Omega)}
    \Big( \sum_{\F \in \FcalnIN}
    \sum_{\TEF \in \tautildeEF}{}
    \Norm{\EfrakEF\jump{v}}_{\W^{\frac{1}{\qprime},\q}
    (\partial\TEF)}^{\q}
    \Big)^\frac{1}{\q}\\
&   \overset{\eqref{eq:stability-zero-extension}}{\le}
    \CtildeTR(\q,\gamma)
    \CHt(\qprime,\Omega,\GammaD)
    \CEG(\q,\gamma)
    \Norm{\sigmabold}_{\Lbf^{\qprime}(\Omega)}
    \Big( \sum_{\F \in \FcalnIN}
    \sum_{\TEF \in \tautildeEF}{}
    \Norm{\jump{v}}_{\W^{\frac{1}{\qprime},\q}
    (\FTEF)}^{\q}
    \Big)^\frac{1}{\q}\\
&   \le \CtildeTR(\q,\gamma)
    \CHt(\qprime,\Omega,\GammaD)
    \CEG(\q,\gamma)
    \Norm{\sigmabold}_{\Lbf^{\qprime}(\Omega)}
    \Big( \sum_{\F \in \FcalnIN} {}
    \Norm{\jump{v}}_{\W^{\frac{1}{\qprime},\q} (\F)}^{\q}
    \Big)^\frac{1}{\q}
\end{split}
\end{equation}\normalsize
Collecting~\eqref{eq:initial-step-W2-broken},
\eqref{eq:initial-split-Helmholtz},
\eqref{eq:T1-version-1&2}, and~\eqref{eq:T2-version-12},
\eqref{eq:sharp-W2p-inequalities-2} follows.

\subsection{Proof of the nonsharp inequalities} \label{subsection:nonsharp-inequalities}
Here, we prove Theorem~\ref{theorem:nonsharp-ineq}.

\paragraph*{A duality rewriting.}
Let~$v$ be in~$\W^{2,\p}(\taun)$, $1 < \p < \infty$.
By duality, we have
\begin{equation} \label{eq:initial-step-W2-broken-nonsharp}
\Norm{\nablaboldh v}_{\Lbf^{\postar}(\Omega)}
= \sup_{\sigmabold \in \Lbf^{\postarprime}(\Omega)}
    \frac{_{\Lbf^{\postarprime}} ( \sigmabold, \nablaboldh v )_{\Lbf^{\postar}}}
        {\Norm{\sigmabold}_{\Lbf^{\postarprime}(\Omega)}}.
\end{equation}
Using decomposition~\eqref{eq:Helmholtz-decomposition} with~$\q=\postar$
we get
\begin{equation} \label{eq:initial-split-Helmholtz-nonsharp}
_{\Lbf^{\postarprime}} ( \sigmabold, \nablaboldh v )_{\Lbf^{\postar}}
    = _{\Lbf^{\postarprime}} ( \nablabold w, \nablaboldh v )_{\Lbf^{\postar}}
    + _{\Lbf^{\postarprime}} ( \xibold, \nablaboldh v )_{\Lbf^{\postar}}
    =: T_1 + T_2.
\end{equation}
We focus on the two terms on the right-hand side separately.

\paragraph*{The term~$T_1$.}
Integrating by parts, we have
\[
T_1
= _{\L^{\pprime}(\Omega)} (w, \Deltah v )_{\L^{\p}(\Omega)}
    + \sum_{\F \in \FcalnID} {}
        _{\L^{\pprime}(\F)} ( w , 
                \jump{\nablaboldh v}_N )_{\L^{\p}(\F)}
=:  T_{1,1} + T_{1,2}.
\]
As for the term~$\T_{1,1}$,
noting that $\pprime\ostar=(\postarprime)^*$,
we preliminarily write
\small\[
\begin{split}
\Norm{w}_{\L^{\pprime}(\Omega)}
&   \le \vert\Omega\vert^{\frac{d-1}{d}}
    \Norm{w}_{\L^{\pprime\ostar}(\Omega)}
    =   \vert\Omega\vert^{\frac{d-1}{d}}
    \Norm{w}_{\L^{(\postarprime)^*}(\Omega)}\\
&   \overset{\eqref{eq:standard-Poincare-Steklov},\eqref{eq:GNS-inequality}}{\le}
    \vert\Omega\vert^{\frac{d-1}{d}}\CSob(\postarprime,\Omega)
    (1+\CPS(\postarprime))
    \SemiNorm{w}_{\W^{1,\postarprime}(\Omega)}.
\end{split}
\]\normalsize
Overall, we write
\[
\begin{split}
T_{1,1} \le \vert\Omega\vert^{\frac{d-1}{d}}
    \CSob(\postarprime,\Omega)
    (1+\CPS(\postarprime))
\SemiNorm{w}_{\W^{1,\postarprime}(\Omega)}
    \Norm{\Deltah v}_{\L^{\p}(\Omega)}.
\end{split}
\]
As for the term~$T_{1,2}$,
using two H\"older's inequalities
with indices $\p$ and $\pprime$,
and~$\p$ and~$\pprime$ (for sequences), we get
\begin{equation} \label{eq:2nd-2nd-term-a}
\begin{aligned}
T_{1,2}
&   =\sum_{\F \in \FcalnID}{}
    _{\L^{\p}(\F)}(\jump{\nablaboldh v}_N,w)_{\L^{\pprime}(\F)}
    \le \sum_{\F \in \FcalnID}{}
        \htildeF^{\frac{1-\p}{\p}} 
        \Norm{\jump{\nablaboldh v}_N}_{\L^{\p}(\F)}
        \htildeF^{-\frac{1-\p}{\p}} \Norm{w}_{\L^{\pprime}(\F)}\\
&   \le \Big( \sum_{\F \in \FcalnID}{} 
        \htildeF^{1-\p}
        \Norm{\jump{\nablaboldh v}_N}_{\L^{\p}(\F)}^{\p} \Big)^\frac{1}{\p}
    \Big( \sum_{\F \in \FcalnID}{} 
            \htildeF \Norm{w}_{\L^{\pprime}(\F)}^{\pprime}    
                        \Big)^\frac{1}{\pprime} .
\end{aligned}
\end{equation}
We focus on the second term on the right-hand side,
since the first one appears in the final estimates.
Using that $\htildeF$ is smaller than or equal to~$\hE$
for $\E$ such that $\F$ belongs to~$\FcalE$ yields
\begin{equation} \label{eq:zuzzu}
\begin{aligned}
&   \Big( \sum_{\F \in \FcalnID}{} 
            \htildeF \Norm{w}_{\L^{\pprime}(\F)}^{\pprime}    
                        \Big)^\frac{1}{\pprime}
    \le \Big( \sum_{\E \in \taun}{}  \hE
            \sum_{\F\in\FcalE}
            \Norm{w}_{\L^{\pprime}(\F)}^{\pprime}    
                        \Big)^\frac{1}{\pprime}
    = \Big(  \hspace{-.5mm}\sum_{\E \in \taun}{}   \hspace{-.5mm}\hE
            \Norm{w}_{\L^{\pprime}(\partial\E)}^{\pprime}
            \Big)^\frac{1}{\pprime} .
\end{aligned}
\end{equation}
The multiplicative trace inequality,
cf. \cite[proof of Theorem 1.2]{Botti-Mascotto:2026},
\begin{equation} \label{eq:multiplicative-trace}
\Norm{w}_{\L^{\q}(\partial \E)}^{\q}
\le  \frac{d}{\gamma\hE}  |\E|^{\frac{s-\q}{s}} 
        \Norm{w}_{\L^{s}(\E)}^{\q} 
      + \frac{\q}{\gamma} 
        \SemiNorm{w}_{\W^{1,\frac{s}{s-\q+1}}(\E)}
        \Norm{w}_{\L^{s}(\E)}^{\q-1},
        \quad\q\in[1,\infty),\;
        \s\in[q,\infty),
\end{equation}
with parameters~$\q=\pprime$ and~$\s=\pprime\ostar$ implies
\[
\begin{aligned}
\hE\Norm{w}_{\L^{\pprime}(\partial\E)}^{\pprime}
&   \le \hE\Norm{w}_{\L^{\pprime}(\partial\E)}^{\pprime}\\
&   \le \frac{d}{\gamma} \vert\E\vert^\frac1d 
            \Norm{w}^{\pprime}_{\L^{\pprime\ostar}(\E)}
        + \frac{\pprime}{\gamma}
            \SemiNorm{w}_{\W^{1,(\p\ostar)'}(\E)}
            \Norm{w}^{\pprime-1}_{\Lbf^{\pprime\ostar}(\E)}.
\end{aligned}
\]
Summing over all elements yields
\[
\begin{aligned}
&   \sum_{\E \in \taun}{} \hspace{-.5mm} \hE
    \Norm{w}_{\L^{\pprime}(\partial\E)}^{\pprime}\\
&   \le \frac{d}{\gamma} \sum_{\E \in \taun}{}   \hspace{-.5mm}
     \vert\E\vert^\frac1d 
            \Norm{w}^{\pprime}_{\L^{\pprime\ostar}(\E)}
    +  \frac{\pprime}{\gamma} \sum_{\E \in \taun}{}   \hspace{-.5mm}
        \SemiNorm{w}_{\W^{1,(\p\ostar)'}(\E)}
        \Norm{w}^{\pprime-1}_{\L^{\pprime\ostar}(\E)}
        =: T_{1,2,1}+T_{1,2,2}.
\end{aligned}
\]
Using H\"older's inequality for sequences
with indices $d$ and $d'=\ostar$ entails
\[
\begin{aligned}
T_{1,2,1}
&   \le \frac{d}{\gamma} 
    \Big( \sum_{\E \in \taun}{} \hspace{-.5mm} 
        \vert\E\vert \Big)^\frac1d
    \Big( \sum_{\E \in \taun}{} \hspace{-.5mm} 
        \Norm{w}_{\L^{\pprime\ostar}(\E)}
            ^{\pprime\ostar} \Big)^\frac1{\ostar}
    = \frac d\gamma \vert\Omega\vert^{\frac1d}
        \Norm{w}^{\pprime}_{\L^{\pprime\ostar}(\Omega)}\\
&   \overset{\eqref{eq:standard-Poincare-Steklov},\eqref{eq:GNS-inequality}}{\le}
    \frac d\gamma \vert\Omega\vert^{\frac1d}
    \big(\CSob(\postarprime,\Omega)
    (1 + \CPS(\postarprime))\big)^{\pprime}
    \SemiNorm{w}_{\W^{1,(\postar)'}(\Omega)}^{\pprime}.
\end{aligned}
\]
Using H\"older's inequality for sequences
with indices $\postarprime$ and $\p\ostar$,
and noting that $(\pprime-1)\postar=\pprimeostar$,
we also have
\[
\begin{aligned}
T_{1,2,2}
&   \le \frac{\pprime}{\gamma} 
    \Big( \sum_{\E \in \taun}{} \hspace{-.5mm}
    \SemiNorm{w}_{\W^{1,\postarprime}(\E)}^{\postarprime} \Big)^\frac1{\postarprime}
    \Big( \sum_{\E \in \taun}{} \hspace{-.5mm} 
        \Norm{w}_{\L^{\pprime\ostar}(\E)}
            ^{\pprime\ostar} \Big)^\frac1{\postar}
    =   \frac{\pprime}{\gamma}
        \SemiNorm{w}_{\W^{1,\postarprime}(\Omega)}
    \Norm{w}_{\L^{\pprime\ostar}(\Omega)}
            ^{\frac{1}{\p-1}}\\
&   \overset{\eqref{eq:standard-Poincare-Steklov},\eqref{eq:GNS-inequality}}{\le}
    \frac{\pprime}{\gamma}
    \big(\CSob(\postarprime,\Omega)
    (1 + \CPS(\postarprime))\big)^{\frac{1}{\p-1}}
        \SemiNorm{w}_{\W^{1,\postarprime}(\Omega)}^{\pprime}.
\end{aligned}
\]
Combining the three displays above implies
\[
\sum_{\E \in \taun}{} \hspace{-.5mm} \hE
\Norm{w}_{\L^{\pprime}(\partial\E)}^{\pprime}
\le \text{C}_{\text{MT}}^{\pprime}(\p,\Omega)
    \SemiNorm{w}_{\W^{1,\postarprime}(\Omega)}^{\pprime},
\]
where
\begin{equation} \label{eq:CMT-definition}
\begin{split}
\text{C}_{\text{MT}}^{\pprime}(\p,\Omega)
&   := \frac d\gamma \vert\Omega\vert^{\frac1d}
    \big(\CSob(\postarprime,\Omega)
    (1 + \CPS(\postarprime))\big)^{\pprime}\\
&   \quad    + \frac{\pprime}{\gamma}
    \big(\CSob(\postarprime,\Omega)
    (1 + \CPS(\postarprime))\big)^{\frac{1}{\p-1}}.
\end{split}
\end{equation}
The bound above and~\eqref{eq:zuzzu} yield
\[
\Big( \sum_{\F \in \FcalnID}{} 
        \htildeF \Norm{w}_{\L^{\pprime}(\F)}^{\pprime}    
        \Big)^\frac{1}{\pprime}
\le \CMT(\p,\Omega) \SemiNorm{w}_{\W^{1,\postarprime}(\Omega)}.
\]
A combination of this with~\eqref{eq:2nd-2nd-term-a} gives
\[
T_{1,2}
\le \CMT(\p,\Omega)
    \SemiNorm{w}_{\W^{1,\postarprime}(\Omega)}
    \Big( \sum_{\F \in \FcalnID}{}
            \htildeF^{1-\p} \Norm{\jump{\nablaboldh v}_N}_{\L^{\p}(\F)}^{\p} \Big)^\frac{1}{\p}.
\]
We combine the bounds on~$T_{1,1}$ and~$T_{1,2}$, and write
\small\begin{equation} \label{eq:T1-nonsharp}
\begin{aligned}
&   T_1
    \le
    \SemiNorm{w}_{\W^{1,\postarprime}(\Omega)}
        \Big[ \vert\Omega\vert^{\frac{d-1}{d}}
        \CSob(\postarprime,\Omega)
        (1+\CPS(\postarprime))
        \Norm{\Deltah v}_{\L^{\p}(\Omega)}\\
&   \qquad\qquad\qquad\qquad\qquad
    + \CMT(\p,\Omega) \Big( \sum_{\F \in \FcalnID}{}
            \htildeF^{1-\p} \Norm{\jump{\nablaboldh v}_N}_{\L^{\p}(\F)}^{\p} \Big)^\frac{1}{\p} \Big]\\
&   \overset{\eqref{eq:regularity-shift-Helmholtz}}{\le}
    \CHo(\postarprime,\Omega,\GammaD) 
    \Norm{\sigma}_{\Lbf^{\postarprime}(\Omega)}
    \Big[ \vert\Omega\vert^{\frac{d-1}{d}}
    \CSob(\postarprime,\Omega)
    (1+\CPS(\postarprime))
    \Norm{\Deltah v}_{\L^{\p}(\Omega)}\\
&   \qquad\qquad\qquad\qquad\qquad
    + \CMT(\p,\Omega) \Big( \sum_{\F \in \FcalnID}{}
            \htildeF^{1-\p} \Norm{\jump{\nablaboldh v}_N}_{\L^{\p}(\F)}^{\p} \Big)^\frac{1}{\p} \Big] .
\end{aligned}
\end{equation}\normalsize

\paragraph*{The term~$T_2$.}
Since the domain~$\Omega$ satisfies~\eqref{eq:assumption-Omega-3},
we consider the regularized right-inverse~$\Psibold$
of the $\curlbold$ operator
as in~\eqref{eq:stability:RI-curl}
of~$\xibold$ and write
\small\[
T_2
    = _{\Lbf^{\postarprime}} (\xibold, \nablaboldh v )_{\Lbf^{\postar}}
    =  \sum_{\E\in\taun}{}
        _{\Lbf^{\postarprime}(\E)} (\curlbold \Psibold,
        \nablaboldh v)_{\Lbf^{\postar}(\E)}
    = \sum_{\E\in\taun}{}
        _{\Lbf^{\pprime}(\partial\E)} (\Psibold,
            \nbfE \times \nablaboldh v)_{\Lbf^{\p}(\partial\E)}.
\]\normalsize
The inner product on the right-hand side
is well posed since $\pprime= (\postarprime)^\sharp$.
We deduce
\[
\begin{split}
T_2
&   = \sum_{\F \in \FcalnIN} {}
        _{\Lbf^{\pprime}(\F)} (\Psibold  ,
            \jump{\nablaboldh v}_T)_{\Lbf^{\p}(\F)} .
\end{split}
\]
Using two H\"older's inequalities
with indices~$\p$ and~$\pprime$,
and~$\p$ and~$\pprime$ (for sequences), we get
\begin{equation} \label{eq:2nd-2nd-term-a-bis}
\begin{aligned}
T_{2}
&   =\sum_{\F \in \FcalnIN}{}
    _{\Lbf^{\p}(\F)}(\jump{\nablaboldh v}_T,\Psibold)_{\Lbf^{\pprime}(\F)}
    \le \sum_{\F \in \FcalnIN}{}
        \htildeF^{\frac{1-\p}{\p}} 
        \Norm{\jump{\nablaboldh v}_T}_{\Lbf^{\p}(\F)}
        \htildeF^{-\frac{1-\p}{\p}} \Norm{\Psibold}_{\Lbf^{\pprime}(\F)}\\
&   \le \Big( \sum_{\F \in \FcalnIN}{} 
        \htildeF^{1-\p}
        \Norm{\jump{\nablaboldh v}_T}_{\Lbf^{\p}(\F)}^{\p} \Big)^\frac{1}{\p}
    \Big( \sum_{\F \in \FcalnIN}{} 
            \htildeF \Norm{\Psibold}_{\Lbf^{\pprime}(\F)}^{\pprime}    
                        \Big)^\frac{1}{\pprime} .
\end{aligned}
\end{equation}
We focus on the second term on the right-hand side,
since the first one appears in the final estimates.
Using that~$\htildeF$ is smaller than or equal to~$\hE$
for $\E$ such that $\F$ belongs to~$\FcalE$ yields
\begin{equation} \label{eq:zuzzu-bis}
\begin{aligned}
&   \Big( \sum_{\F \in \FcalnIN}{} 
            \htildeF \Norm{\Psibold}_{\Lbf^{\pprime}(\F)}^{\pprime}    
                        \Big)^\frac{1}{\pprime}
    \le \Big( \sum_{\E \in \taun}{}  \hE
            \sum_{\F\in\FcalE}
            \Norm{\Psibold}_{\Lbf^{\pprime}(\F)}^{\pprime}    
                        \Big)^\frac{1}{\pprime}
    = \Big(  \hspace{-.5mm}\sum_{\E \in \taun}{}   \hspace{-.5mm}\hE
            \Norm{\Psibold}_{\Lbf^{\pprime}(\partial\E)}^{\pprime}
            \Big)^\frac{1}{\pprime}.
\end{aligned}
\end{equation}
Using the multiplicative trace
inequality~\eqref{eq:multiplicative-trace} with
parameters~$\q=\pprime$ and~$\s=\pprime\ostar$ gives
\[
\begin{aligned}
\hE\Norm{\Psibold}_{\Lbf^{\pprime}(\partial\E)}^{\pprime}
&   \le \hE\Norm{\Psibold}_{\Lbf^{\pprime}(\partial\E)}^{\pprime}\\
&   \le \frac{d}{\gamma} \vert\E\vert^\frac1d 
            \Norm{\Psibold}^{\pprime}_{\Lbf^{\pprime\ostar}(\E)}
        + \frac{\pprime}{\gamma}
            \SemiNorm{\Psibold}_{\Wbf^{1,(\p\ostar)'}(\E)}
            \Norm{\Psibold}^{\pprime-1}_{\Lbf^{\pprime\ostar}(\E)}.
\end{aligned}
\]
Summing over all elements this estimate
further yields
\[
\begin{aligned}
&   \sum_{\E \in \taun}{} \hspace{-.5mm} \hE
    \Norm{\Psibold}_{\Lbf^{\pprime}(\partial\E)}^{\pprime}\\
&   \le \frac{d}{\gamma} \sum_{\E \in \taun}{}   \hspace{-.5mm}
     \vert\E\vert^\frac1d 
            \Norm{\Psibold}^{\pprime}_{\Lbf^{\pprime\ostar}(\E)}
    +  \frac{\pprime}{\gamma} \sum_{\E \in \taun}{}   \hspace{-.5mm}
        \SemiNorm{\Psibold}_{\Wbf^{1,(\p\ostar)'}(\E)}
        \Norm{\Psibold}^{\pprime-1}_{\Lbf^{\pprime\ostar}(\E)}
        =: T_{2,1}+T_{2,2}.
\end{aligned}
\]
Using H\"older's inequality for sequences
with indices $d$ and $d'=\ostar$
entails
\[
\begin{aligned}
T_{2,1}
&   \le \frac{d}{\gamma} 
    \Big( \sum_{\E \in \taun}{} \hspace{-.5mm} 
        \vert\E\vert \Big)^\frac1d
    \Big( \sum_{\E \in \taun}{} \hspace{-.5mm} 
        \Norm{\Psibold}_{\Lbf^{\pprime\ostar}(\E)}
            ^{\pprime\ostar} \Big)^\frac1{\ostar}
    = \frac d\gamma \vert\Omega\vert^{\frac1d}
        \Norm{\Psibold}^{\pprime}_{\Lbf^{\pprime\ostar}(\Omega)}\\
&   \le \frac d\gamma \vert\Omega\vert^{\frac1d}
        \big(\CSob(\postarprime,\Omega)
        (1 + \CPS(\postarprime))\big)^{\pprime} 
        \SemiNorm{\Psibold}_{\Wbf^{1,(\postar)'}(\Omega)}^{\pprime}.
\end{aligned}
\]
Using H\"older's inequality for sequences
with indices $\postarprime$ and $\p\ostar$,
and noting that $(\pprime-1)\postar=\pprimeostar$,
we also have
\[
\begin{aligned}
T_{2,2}
&   \le \frac{\pprime}{\gamma} 
    \Big( \sum_{\E \in \taun}{} \hspace{-.5mm}
        \SemiNorm{\Psibold}_{\Wbf^{1,\postarprime}(\E)}^{\postarprime} \Big)^\frac1{\postarprime}
    \Big( \sum_{\E \in \taun}{} \hspace{-.5mm} 
        \Norm{\Psibold}_{\Lbf^{\pprime\ostar}(\E)}
            ^{\pprime\ostar} \Big)^\frac1{\postar}
    \le \frac{\pprime}{\gamma}
        \SemiNorm{\Psibold}_{\Wbf^{1,\postarprime}(\E)}
    \Norm{\Psibold}_{\Lbf^{\pprime\ostar}(\Omega)}
            ^{\frac{1}{\p-1}}\\
&   \overset{\eqref{eq:standard-Poincare-Steklov},\eqref{eq:GNS-inequality}}{\le}
    \frac{\pprime}{\gamma}
    \big(\CSob(\postarprime,\Omega)
    (1 + \CPS(\postarprime))\big)^{\frac{1}{\p-1}}
        \SemiNorm{\Psibold}_{\Wbf^{1,\postarprime}(\Omega)}^{\pprime}.
\end{aligned}
\]
Given~$\CMT$ as in~\eqref{eq:CMT-definition},
combining the three displays above implies
\[
\sum_{\E \in \taun}{} \hspace{-.5mm} \hE
\Norm{\Psibold}_{\Lbf^{\pprime}(\partial\E)}^{\pprime}
\le \text{C}_{\text{MT}}^{\pprime}(\p,\Omega)
    \SemiNorm{\Psibold}_{\Wbf^{1,\postarprime}(\Omega)}^{\pprime}.
\]
The bound above and~\eqref{eq:zuzzu-bis} yield
\[
\Big( \sum_{\F \in \FcalnIN}{} 
        \htildeF \Norm{w}_{\L^{\pprime}(\F)}^{\pprime}    
        \Big)^\frac{1}{\pprime}
\le \CMT(\p,\Omega) \SemiNorm{w}_{\W^{1,\postarprime}(\Omega)}
\]
A combination of this display with~\eqref{eq:2nd-2nd-term-a-bis} gives
\footnotesize\begin{equation} \label{eq:T2-nonsharp}
\begin{split}
T_{2}
&   \le \sum_{\F \in \FcalnIN}{}
    _{\L^{\p}(\F)}(\jump{\nablaboldh v}_T, w)_{\L^{\pprime}(\F)}
    \le \CMT(\p,\Omega)
    \SemiNorm{\Psibold}_{\Wbf^{1,\postarprime}(\Omega)}
    \Big( \sum_{\F \in \FcalnIN}{}
            \htildeF^{1-\p}
            \Norm{\jump{\nablaboldh v}_T}
            _{\Lbf^{\p}(\F)}^{\p} \Big)^\frac{1}{\p}\\
&   \overset{\eqref{eq:stability:RI-curl}}{\le}
    \CMT(\p,\Omega)
    \CBcal(\postarprime,\Omega,\GammaD)
    \CHt(\postarprime,\Omega,\GammaD)
    \Norm{\sigmabold}_{\Lbf^{\postarprime}(\Omega)}
    \Big( \sum_{\F \in \FcalnIN}{}
            \htildeF^{1-\p} \Norm{\jump{\nablaboldh v}_T}_{\Lbf^{\p}(\F)}^{\p} \Big)^\frac{1}{\p}.
\end{split}
\end{equation}\normalsize

\paragraph*{Wrap up.}
Combining~\eqref{eq:initial-step-W2-broken-nonsharp},
\eqref{eq:initial-split-Helmholtz-nonsharp},
\eqref{eq:T1-nonsharp}, and~\eqref{eq:T2-nonsharp},
we arrive at
\small\[
\begin{split}
\Norm{\nablaboldh v}_{\Lbf^{\postar}(\Omega)}
&   \le \CHo(\postarprime,\Omega,\GammaD)
    \Big[ \vert\Omega\vert^{\frac{d-1}{d}}
    \CSob(\postarprime,\Omega)
    (1+\CPS(\postarprime))
        \Norm{\Deltah v}_{\L^{\p}(\Omega)}\\
&   \quad +  \CMT(\p,\Omega)
            \Big( \sum_{\F \in \FcalnID}{}
            \htildeF^{1-\p} \Norm{\jump{\nablaboldh v}_N}_{\L^{\p}(\F)}^{\p} \Big)^\frac{1}{\p} \Big] \\
&   \quad + \CMT(\p,\Omega)
            \CBcal(\postarprime,\Omega,\GammaD)
            \CHt(\postarprime,\Omega,\GammaD)
            \Big( \sum_{\F \in \FcalnIN}{}
            \htildeF^{1-\p} \Norm{\jump{\nablaboldh v}_T}_{\Lbf^{\p}(\F)}^{\p} \Big)^\frac{1}{\p},
\end{split}
\]\normalsize
which is~\eqref{eq:nonsharp-W2p-inequalities-1}.

\paragraph*{The term~$T_2$: alternative approach.}
We now assume that~$\Omega$ has general topology
and $\p\in(1,2(d-1)/d)$.
Proceeding as in the derivation of~\eqref{eq:sharp-W2p-inequalities-2}
taking $\q=\postar$, \eqref{eq:nonsharp-W2p-inequalities-2} follows.

\paragraph*{Acknowledgments.}
The authors are members of the Gruppo Nazionale Calcolo Scientifico-Istituto Nazionale di Alta Matematica (GNCS-INdAM).

\paragraph*{Fundings.}
The authors have been partially funded by the
European Union (ERC, NEMESIS, project number 101115663);
views and opinions expressed are however those
of the authors only and do not necessarily reflect
those of the EU or the ERC Executive Agency.

\paragraph*{Declaration.}
\begin{quote}
\itshape
``Instinct, however well honed, could not necessarily warn one
that a young killer, hidden behind a rock four hundred yards away\dots\
was looking through a telescope sight, about to squeeze the trigger.''

\smallskip
\raggedleft \textbf{Larry McMurtry}, \emph{Streets of Laredo}
\end{quote}

\noindent
A final exercise in traditional craftsmanship before the sight takes over.

{\footnotesize
\bibliography{bibliography.bib}}
\bibliographystyle{plainurl}

\appendix

\section{Geometric singularities of solutions to the Poisson problem} \label{appendix:regularity-Poisson}
In the analysis of Section~\ref{section:main},
it will be crucial to know the $\W^{1,\q}$-regularity of solutions
to problem~\eqref{eq:1st-element-Helmholtz},
which depends on the choice of boundary conditions.
Below, we review relevant regularity results from the literature.

\paragraph*{Lipschitz domains with full Dirichlet or full Neumann boundary conditions.}
In spatial dimensions $d \in \{2,3\}$, for sufficiently smooth data,
the solution to~\eqref{eq:1st-element-Helmholtz}
with full Neumann boundary conditions ($\Gamma = \Gamma_{\mathrm{N}}$)
belongs to $\W^{1,\q}(\Omega)$ for
\begin{equation} \label{eq:JK-bounds}
\q \in \left( \frac{2d}{d+1}, \frac{2d}{d-1} \right);
\end{equation}
see, e.g., \cite{Geng:2012, Jerison-Kenig:1981}.
For full Dirichlet boundary conditions, the same interval of~$\q$ holds;
see, e.g., \cite{Jerison-Kenig:1995} and~\cite{Shen:2005}.
In fact, these results hold for
$\q \in \big(\frac{2d}{d+1}-\varepsilon, \frac{2d}{d-1}+\varepsilon\big)$
for a domain-dependent $\varepsilon > 0$.
These $\W^{1,\q}(\Omega)$ bounds are in direct agreement
with the $H^{\frac{3}{2}}(\Omega)$ regularity established
in \cite[Theorems~3 and~4]{Savare:1998} via Sobolev embeddings.

\begin{remark} \label{remark:reg-mixed-BCs}
The case of mixed boundary conditions is less addressed in the literature;
we refer to \cite{Groeger:1989, HallerDintelmann-Kaiser-Rehberg:2008,
Taylor-Ott-Brown:2013, Ott-Brown:2013} for results in this direction.
Because boundary singularities arise at the interface
where boundary conditions change,
the admissible range of~$\q$ is generally
narrower than in~\eqref{eq:JK-bounds}.
\end{remark}

\paragraph*{Polygonal domains with either full Dirichlet or full Neumann boundary conditions.}
In this case (see, e.g., \cite[Theorem~4.1]{Schwab:1998}),
the solution to~\eqref{eq:1st-element-Helmholtz}
can be decomposed into a smooth part (related to the regularity shift)
and a singular part consisting of linear combinations
of corner singularities.
The leading singular term has the form
$S(r,\theta)=r^\alpha \varphi(\theta)$,
where $\varphi(\theta)$ is a smooth function
in the angular variable~$\theta$ and $\alpha=\pi/\omega$,
with~$\omega$ being the largest interior angle of~$\Omega$.
If $\alpha$ is an integer, the singular radial part
reads $r^\alpha \log(r)$ instead of~$r^\alpha$.
Direct computations reveal that
\[
S(r,\theta) \in
\begin{cases}
\H^{1+\alpha-\varepsilon}(\Omega)
    \cap \W^{1,\frac{2}{1-\alpha}-\varepsilon}(\Omega)
    & \text{if } \omega > \pi,\\[1ex]
\H^{1+\alpha-\varepsilon}(\Omega)
    \cap \W^{1,\infty}(\Omega)
    & \text{if } \omega \le \pi,
\end{cases}
\qquad\qquad\forall\varepsilon>0.
\]
We outline the $\W^{1,\q}$-regularity of~$S$. First, observe that
\[
\partial_r S(r,\theta) = \alpha r^{\alpha-1} \varphi(\theta),
\qquad\qquad
r^{-1}\partial_\theta S(r,\theta) = r^{\alpha-1} \varphi'(\theta),
\]
whence
\[
\vert \nablabold S \vert
= \sqrt{(\partial_r S(r,\theta))^2 + (r^{-1}\partial_\theta S(r,\theta))^2}
= r^{\alpha-1} \sqrt{\alpha^2 \varphi(\theta)^2 + \varphi'(\theta)^2}.
\]
Focusing on a circular neighbourhood~$B_\omega$ of radius~$R$
centered at the vertex with the largest angle~$\omega$, we deduce
\begin{equation} \label{eq:computations-W1q-regularity}
\SemiNorm{S}_{\W^{1,\q}(\Omega\cap B_\omega)}^{\q}
= \Big( \int_{0}^\omega \big(
    \alpha^2\varphi(\theta)^2 + \varphi'(\theta)^2
        \big)^\frac\q2 d\theta \Big) 
        \int_{0}^R r^{(\alpha-1)\q +1} dr
    = C_\omega
        \int_{0}^R r^{(\alpha-1)\q +1} dr.
\end{equation}
The integral on the right-hand side is finite for $\alpha<1$
if $(\alpha-1)\q+1>-1$, i.e., $\q < \frac{2}{1-\alpha}$.
Of course, the lower bound in~\eqref{eq:JK-bounds} must hold as well.

\paragraph*{Polygonal domains with mixed Dirichlet--Neumann boundary conditions.}
Following again \cite[Theorem~4.1]{Schwab:1998},
similar comments apply to the mixed boundary condition case.
The main distinction is that the singularity exponent
is now $\alpha = \frac{\pi}{2\omega}$,
where~$\omega$ is the largest interior angle of~$\Omega$,
specifically at an interface where boundary conditions change.
The leading singular term has either the form $r^{\alpha}\varphi(\theta)$
or $r^{\alpha} \log(r)\varphi(\theta)$,
depending on whether~$\alpha$ is an integer.
Direct computations reveal that
\[
S(r,\theta) \in
\begin{cases}
\H^{1+\alpha-\varepsilon}(\Omega)
    \cap \W^{1,\frac{2}{1-\alpha}-\varepsilon}(\Omega)
    & \text{if } \omega > \frac{\pi}{2},\\[1ex]
\H^{1+\alpha-\varepsilon}(\Omega)
    \cap \W^{1,\infty}(\Omega)
    & \text{if } \omega \le \frac{\pi}{2},
\end{cases}
\qquad\qquad\qquad\forall\varepsilon>0.
\]
This regularity can be derived via computations
analogous to those in~\eqref{eq:computations-W1q-regularity}.
A lower bound on the regularity exponent must also be taken into account;
cf. Remark~\ref{remark:reg-mixed-BCs}.

\paragraph*{Polyhedral domains with either full Dirichlet or full Neumann boundary conditions.}
The types of singularities in polyhedra are more involved
than in the polygonal case; see, e.g., \cite{Dauge:1988}.
Here, there are essentially three types of singularities:
\begin{itemize}
\item \emph{Vertex singularities}
of the form $r_V^{\lambda_V} \psi(\theta_V,\varphi_V)$,
where $\psi(\cdot,\cdot)$ is a smooth function
in the two angular components $\theta_V$ and $\varphi_V$,
$r_V$ is the radial variable associated with the vertex~$V$,
and $\lambda_V > 0$ is related to the smallest positive
(Dirichlet or Neumann) eigenvalue of the Laplace-Beltrami operator
on the intersection of the local tangent cone
departing from the vertex with the unit sphere;
\item \emph{Edge singularities}, where the singular behaviour
occurs in the direction orthogonal to the edge and matches
the vertex singularities of 2D polygons, i.e.,
they have the form $r_e^{\alpha_e} \chi(\theta_e) g_e(z_e)$
with $\alpha_e = \frac{\pi}{\omega_e}$,
where $\omega_e$ is the opening angle at the edge~$e$
(with the usual logarithmic factor for integer $\alpha_e$),
$(r_e,\theta_e)$ denote the 2D polar coordinates
in the plane orthogonal to the edge~$e$,
and $g_e(z_e)$ is a smooth factor along the edge coordinate~$z_e$;
\item \emph{Vertex--edge singularities},
which occur at the transition where an edge meets a vertex;
in particular, the edge factor $g_e(z_e)$ develops
a singular behaviour when approaching
the vertices~$V_j$ ($j=1,2$) of the edge~$e$.
For each edge, there are two extra interaction singularities of the form
\[
r_{V_j}^{\lambda_{V_j}} 
    \Big( \frac{r_e}{r_{V_j}} \Big)^{\alpha_e}
    \psi(\theta_e).
\]
\end{itemize}
Estimating the $\W^{1,\q}$-regularity in this context
is more challenging and requires knowing
the vertex exponents $\lambda_{V_j}$,
which are generally not available in closed form,
in addition to accounting for the lower bound in~\eqref{eq:JK-bounds}.

\paragraph*{Convex polyhedral domains with either full Dirichlet or full Neumann boundary conditions.}
In this case, it is known that the solution
to~\eqref{eq:1st-element-Helmholtz} belongs to $\H^2(\Omega)$,
and therefore to $\W^{1,\q}(\Omega)$ for all $\q \in [1, \infty)$ if $d=2$,
or for all $\q \in [1, 6]$ if $d=3$ via Sobolev embedding.
This $H^2$ regularity fails if mixed boundary conditions are considered,
due to the presence of Zaremba-type singularities
at the interface where boundary conditions change.
Of course, the lower bound in~\eqref{eq:JK-bounds} must hold as well.

\end{document}